\documentclass{elsarticle}

\biboptions{sort&compress,numbers}
\usepackage[nodots]{numcompress}
\usepackage{lineno,url}
\usepackage{epstopdf}
\usepackage{amsmath,amsfonts,amssymb}
\usepackage{amssymb,amsmath,amsthm}
\usepackage{overpic}
\usepackage{contour}\contourlength{0.7pt}
\usepackage{tikz}\usetikzlibrary{shapes,arrows}
\usepackage{versions}
\usepackage{natbib}
\usepackage{algorithm, algorithmic}
\usepackage{multirow}

\newtheoremstyle{mytheorem}{5pt plus 5pt minus 3pt}{4pt plus 3pt minus 1.5pt}
	{\itshape}{}{\bfseries}{.}{1ex plus 1ex minus .5ex}{}
\newtheoremstyle{mydef}{5pt plus 5pt minus 3pt}{4pt plus 3pt minus 1.5pt}
	{}{0pt}{\bfseries}{.}{1ex plus 1ex minus .5ex}{}
\newtheoremstyle{myremark}{5pt plus 5pt minus 3pt}{4pt plus 3pt minus 1.5pt}
	{}{0pt}{\itshape}{.}{1ex plus 1ex minus .5ex}{}
\theoremstyle{mytheorem}
\newtheorem{prop}{Proposition}[section]

\theoremstyle{mydef}

\newtheorem{ex}[prop]{Example}
\theoremstyle{myremark}
\newtheorem{rem}{Remark}

\newcommand{\bp}{\mathbf p}
\newcommand{\de}{\mathrm{d}}

\newcommand{\D}{\mathbf D}
\newcommand{\DDt}{\widetilde{\mathbf D}}
\newcommand{\Q}{\mathcal Q}
\newcommand{\Qt}{\widetilde{\Q}}
\newcommand{\taut}{\widetilde{\tau}}
\newcommand{\omegat}{\widetilde{\omega}}
\newcommand{\tf}{\footnotesize}

\newcommand{\bx}{\mathbf x}
\newcommand{\br}{\mathbf r}
\def\R {\mathbb R}
\newcommand{\FF}{\mathcal F}
\newcommand{\FFt}{\widetilde{\mathcal F}}

\newcommand{\Om}{\Omega}
\newcommand{\Omt}{\widetilde{\Omega}}
\newcommand{\HH}{\mathcal H}
\newcommand{\XX}{\mathcal X}
\newcommand{\XXt}{\widetilde{\mathcal X}}
\newcommand{\xt}{\widetilde{x}}
\newcommand{\St}{\widetilde{S}}
\newcommand{\Dt}{\widetilde{D}}

\def\Acaption#1#2{\caption{#2}\vspace*{-#1}}

\definecolor{Mgreen}{RGB}{34,139,34}
\definecolor{blau}{rgb}{0.15,0.2,0.5}
\definecolor{gray}{rgb}{0.5,0.5,0.5}
\definecolor{drot}{rgb}{0.7,0,0.1}
\definecolor{gelb}{rgb}{.55,.40,.1}
\definecolor{magenta}{rgb}{1.,0.,1.}
\definecolor{cyan}{rgb}{0.,1.,1.}
\definecolor{green}{rgb}{0.,1.,0.}
\definecolor{Morange}{rgb}{1.,0.5,0.}

\begin{document}

\begin{frontmatter}
\title{Gaussian quadrature for splines via homotopy continuation: rules for $C^2$ cubic splines}

\author[NumPor]{Michael Barto\v{n}\corref{cor1}}
\ead{Michael.Barton@kaust.edu.sa}
\author[NumPor,ErSE]{Victor Manuel Calo}
\ead{Victor.Calo@kaust.edu.sa}

\cortext[cor1]{Corresponding author}

\address[NumPor]{Numerical Porous Media Center,\\ King Abdullah University of Science and Technology, Thuwal 23955-6900, KSA} 
\address[ErSE]{Applied Mathematics $\&$ Computational Science and Earth Science $\&$ Engineering,\\ King Abdullah University of Science and Technology, Thuwal 23955-6900, KSA}

\begin{abstract}
We introduce a new concept for generating optimal quadrature rules for splines.
Given a target spline space where we aim to generate an optimal quadrature rule,
we build an associated source space with known optimal quadrature and
transfer the rule from the source space to the target one, preserving the number of
quadrature points and therefore optimality.
The quadrature nodes and weights are, considered as a higher-dimensional
point, a zero of a particular system of polynomial equations.
As the space is continuously deformed by modifying the source knot vector, the quadrature
rule gets updated using polynomial homotopy continuation.
For example, starting with $C^1$ cubic splines with uniform knot sequences,
we demonstrate the methodology by deriving the optimal rules for uniform $C^2$ cubic spline spaces
where the rule was only conjectured heretofore.
We validate our algorithm by showing that the resulting quadrature rule is independent of the
path chosen between the target and the source knot vectors as well as the source rule chosen.

%
\end{abstract}


\end{frontmatter}

\section{Introduction}\label{intro}

Numerical integration of univariate functions is a fundamental mathematical task which
is a subroutine of many complex algorithms and is typically frequently invoked.
Naturally, such an integration (or quadrature) rule
must be as efficient as possible. 
We derive a new class of quadrature rules 
that are optimal in the sense that they require the minimal number of function's evaluations.

A quadrature rule, or shortly a \emph{quadrature}, is an \textit{$m$-point rule},
if $m$ evaluations of a function $f$ are needed to approximate its weighted integral
over a closed interval $[a,b]$
\begin{equation}\label{eq:GaussQuad}
\int_a^b w(x) f(x) \, \mathrm{d}x = \sum_{i=1}^{m} \omega_i f(\tau_i) + R_{m}(f),
\end{equation}
where $w$ is a fixed non-negative \emph{weight function} defined over $[a,b]$.
The rule is required to be \emph{exact}, that is, $R_m(f) \equiv 0$
for each element of a predefined linear function space $\mathcal{L}$.
The rule is said to be \emph{optimal} if $m$ is the minimal number of
\emph{weights} $\omega_i$ and \textit{nodes} $\tau_i$, points at which $f$ has to be evaluated.

For the space of polynomials, the optimal rule is known to be the classical Gaussian quadrature \cite{Gautschi-1997}
with the order of exactness $2m-1$, that is, only $m$ evaluations are needed to exactly integrate
any polynomial of degree at most $2m-1$. Consider a sequence of polynomials $(p_0,p_1,\ldots,p_m,\ldots )$  that form an orthogonal basis
with respect to the scalar product
\begin{equation}
<f,g> = \int_a^b f(x)g(x)w(x) \de x.
\end{equation}
The quadrature points are the roots of the $m$-th orthogonal polynomial $p_m$ which
in the case when $w(x)\equiv 1$ is the degree-$m$ Legendre polynomial \cite{Szego-1936}.


In this paper, we focus on piece-wise cubic polynomials, cubic splines, but 
the methodology is general and could be used for higher degrees as well.
A univariate space of cubic splines is uniquely determined by its \emph{knot vector}. This knot vector is
a non-decreasing sequence of real numbers called \emph{knots} and the multiplicity of each knot determines the smoothness
between the cubic pieces. To simplify the argument, we first study uniform knot vectors with all interior knots
with \emph{uniform multiplicity}. In particular, we investigate knot vectors with single interior knots which yield $C^2$ cubic
splines and knot vectors with double interior knots which give $C^1$ cubic splines.
For a more detailed introduction on splines,
we refer the reader, e.g., to \cite{deBoor-1972,Hoschek-2002-CAGD,Elber-2001}.

The quadrature rules for splines have been studied since late 50's \cite{Schoenberg-1958,Micchelli-1972,Micchelli-1977}.
Firstly conjectured by Schoenberg \cite{Schoenberg-1958}, later proved by Micchelli and Pinkus \cite{Micchelli-1977},
the conditions of the existence and uniqueness of the optimal (Gaussian) rule have been derived.
For spline spaces with maximum continuity (e.g., $C^2$ cubic splines), Micchelli and Pinkus \cite{Micchelli-1977} proved that there always exists
a \emph{unique} Gaussian quadrature rule.
Their result, however, also reveals a nice phenomenon: the number of optimal nodes stays fixed as long as the number of interior knots
stays constant. Therefore if one desires to derive a class of quadrature rules with the same number of nodes,
spline spaces with higher continuity must have adequately more sub-intervals (elements)
than spaces of lower continuity. This is natural since splines of lower continuity are
limits of the higher continuous ones, when merging continuously two (or several) knots together, and their result is
in agreement with this fact.

For spaces with lower continuity, or when boundary constraints are involved, the rule is not guaranteed to be unique.
Only spaces with fixed continuities were studied in \cite{Micchelli-1977}.
To the best of our knowledge, the results on existence and uniqueness are not known for spaces with mixed continuities,
we refer to these as knot vectors with mixed multiplicities. For example, for cubic splines,
this includes to knot vectors with both single and double knots.
In this work, we derive quadratures for this kind of mixed continuity spaces and show numerically that
optimal quadratures exist.

\begin{figure}[tb]\color{blau!50}
\noindent{\fboxsep.015\columnwidth\fbox{\color{blau}\relax
	\begin{minipage}{.96\columnwidth}\footnotesize
 \tikzstyle{block} = [rectangle, draw, fill=drot!5,
    text width=14em, text centered, rounded corners, minimum height=2em]
 \tikzstyle{block2} = [rectangle, draw, fill=drot!5,
    text width=10em, text centered, rounded corners, minimum height=2em]
 \tikzstyle{info} = [rectangle, color=black, text width=18em, minimum height=2em]
 \tikzstyle{line} = [draw, -latex']

\smallskip
\begin{tikzpicture}[node distance = 1cm]
    \node [block] at (4,-1) (velo)
	{Target space $S$ defined by knot vector $\XX$ for which we derive an optimal quadrature rule};
    \node [block] at (4,-2.5) (knots)
	{$\XX$ and $\XXt$ have the same number of interior knots on $[a,b]$\\ (including multiplicities)};
    \node [block] at (4,-4) (optim)
	{Source space $\St$ defined by knot vector $\XXt$ with a known quadrature rule $\Qt$};
    %
    \path [line] (velo) -- (knots);
    \path [line] (knots) -- (optim);
    \draw[line] (optim) .. controls  (-0,-4) and (-0,-1) .. (velo);
    \draw[line] (velo) .. controls  (8,-1) and (8,-4) .. (optim);
    \node [block2] at (0,-2.5) (homo)
	{Homotopy continuation};


    \node [block2] at (8,-2.5) (initial)
	{Source space initialization};
 \end{tikzpicture}
 \end{minipage}}}\color{black}
 \caption{Spline space initialization flowchart}
 \label{fig:flowchart}
\end{figure}
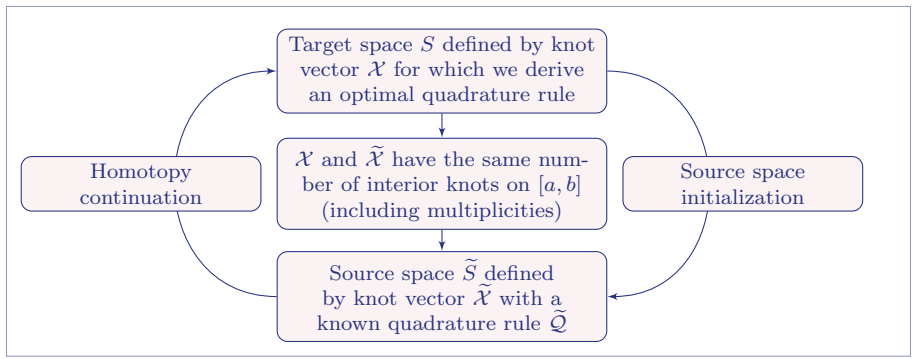

We use the observation that a spline space with multiple knots is a limit case of a spline space with the same number of single knots
(when counting the multiplicities) when two, or several, knots merge. 
Abstracting this merging as a continuous transition between the two knot vectors of the same cardinality,
the \emph{source} and the \emph{target} ones, the corresponding spline spaces continuously evolve from one into the other.
The quadrature rule also depends \emph{continuously} on the spline space
since the quadrature rule can be seen as a zero (root) of a certain polynomial
system and an infinitesimal change of the system does not significantly change the root.
Based on these facts, we propose a new methodology that for a given (target)
spline space $S$ generates an associated source space $\St$
where the Gaussian quadrature rule is known. These spaces are defined above knot vectors $\XX$ and $\XXt$,
respectively, and the quadrature rule $\Qt$ of $\St$ is numerically traced as $\XXt$ evolves into $\XX$,
see Fig.~\ref{fig:flowchart}.

The rest of the paper is organized as follows. Section~\ref{sec:Homotopy} briefly overviews the homotopy continuation of polynomial systems
and Section~\ref{sec:Splines} summarizes a few basic properties of cubic splines. In Section~\ref{sec:HomotopyQuad}, we introduce a
homotopy-continuation-based algorithm and discuss the results and the validity of the new quadratures obtained in Section~\ref{sec:Ex}.
Finally, Section~\ref{sec:conl} summarizes our conclusions and describes future research directions.

\section{Homotopy continuation for polynomial systems}\label{sec:Homotopy}

Polynomial homotopy continuation (PHC) is a numerical scheme that solves polynomial systems of equations.
This approach was introduced to solve the problem of movability of kinematic mechanisms \cite{Wampler-1990} where the variables are the free parameters
of a certain mechanism tied together by a set of polynomial constraints. However, the unknowns and the constraints may relate to an arbitrary problem.
As our Gaussian quadrature rules are derived using a variation of this method, for the sake of completeness,
we now briefly review the ideas used in homotopy continuation.
For a detailed explanation, we refer the reader to the book \cite{Wampler-2005}.

Consider a well-constrained 
$2m \times 2m$ polynomial system
\begin{equation}\label{eq:AlgSys}
     \FF(\bx)=\mathbf 0, \quad \bx \in \Omega \subset \R^{2m},
\end{equation}
where the domain $\Omega$ is a hypercube in $\R^{2m}$, i.e.,
$\Omega =[\underline{x}_1,\overline{x}_1]\times \dots \times [\underline{x}_{2m},\overline{x}_{2m}]$.

Let $K$ be an upper bound of the number of \emph{real} roots of $\FF$ in $\Om$ and
let us consider a  simpler system $\FFt(\bx) = \mathbf 0$ defined over some $\Omt$ with \emph{exactly} $K$ real roots
$\{\br_1,\dots, \br_K\}$, that is,
\begin{equation}\label{eq:AlgSysT}
    \FFt(\br_i) = \mathbf 0, \quad i=1,\dots,K
\end{equation}
For now, let us assume $\Om=\Omt$.
The roots $\{\br_1,\dots, \br_K\}$ are known, in fact they are chosen as an input since 
the system $\FFt$ is as simple as possible, see \cite{Wampler-2005}.
Consider now a continuous deformation of the known system into
the desired one
\begin{equation}\label{eq:ContDef}
    \FFt(\bx)\rightarrow \FF(\bx). 
\end{equation}
Let us denote by $\HH$ a one parameter family of systems created, as $\FFt$ is transformed into $\FF$.
The parameter characterizing the transformation can be thought of as time or pseudo-time $t$. Thus
\begin{equation}\label{eq:HC}
     \HH(\bx,t)=\mathbf 0, \quad \bx \in \Omega, \quad t \in [0,1],
\end{equation}
where $\HH(\bx,0)$ is the artificially constructed system $\FFt(\bx)=\mathbf 0$,
for which the roots are known, and $\HH(\bx,1)$ is the
given system (\ref{eq:AlgSys}) which we seek to solve.
Consider now the roots (\ref{eq:AlgSysT}) as a function of time, $\br_1(t),\dots, \br_K(t)$,
and think of them as trajectories (curves) in $\R^{2m}$. As $t$ runs over $[0,1]$, some of these roots may vanish,
which corresponds to the fact that (\ref{eq:AlgSys}) may have less real roots than $K$,
but no new trajectory may rise. The transformation is guaranteed not to introduce new roots.

Typically, the source and the target systems are blended in a linear fashion, corresponding to the shortest
path when deforming one polynomial system into another.
However, one can consider different paths between the systems \cite{Wampler-2005}.

\section{Transition between $C^1$ and $C^2$ cubic splines}\label{sec:Splines}

We recall several properties of spline basis functions. We
consider a uniform knot vector
\begin{equation}\label{eq:XXt}
\XXt_n = (a=\xt_0,\xt_1,\xt_1,...,\xt_{n-1},\xt_{n-1},\xt_{n} = b)
\end{equation}
on the interval $[a,b]$,
that is, each of the $n-1$ interior knots has multiplicity two.
We define $h := \frac{b-a}{n} = \xt_{k} -  \xt_{k-1}$ for all $k=1, \dots, n$.
Let us denote $N:=2n-1$ and consider a uniform knot vector
\begin{equation}\label{eq:XX}
\XX_N = (a=x_0,x_1,\dots,x_{N-1},x_{N} = b),
\end{equation}
each knot of multiplicity exactly one and define $H := \frac{b-a}{N} = x_{k} -  x_{k-1}$ for all $k=1, \dots, N$.
We denote by $\pi_3$ the space of polynomials of degree at most 3 and define $\St^{n}_{3,1}$,
the linear space of cubic splines over a uniform knot sequence $\XXt_n$ as
\begin{equation}\label{eq:familyC1}
\St^{n}_{3,1} = \{ f\in C^{1}[a,b]: f|_{(\xt_{k-1},\xt_{k})} \in \pi_3, k=1,...,n\}.
\end{equation}
Similarly we denote by $S^{N}_{3,2}$ the linear space of cubic splines over a
uniform knot sequence $\XX_N$
\begin{equation}\label{eq:familyC2}
S^{N}_{3,2} = \{ f\in C^{2}[a,b]: f|_{(x_{k-1},x_{k})} \in \pi_3, k=1,...,N\}.
\end{equation}
The dimension of both spaces is $2n+2=N+3$. That is, the total number of interior knots is the same for both spaces,
while the number of non-zero knot spans is different.

\begin{figure}[!tb]
\vrule width0pt\hfill
 \begin{overpic}[width=.32\textwidth,angle=0]{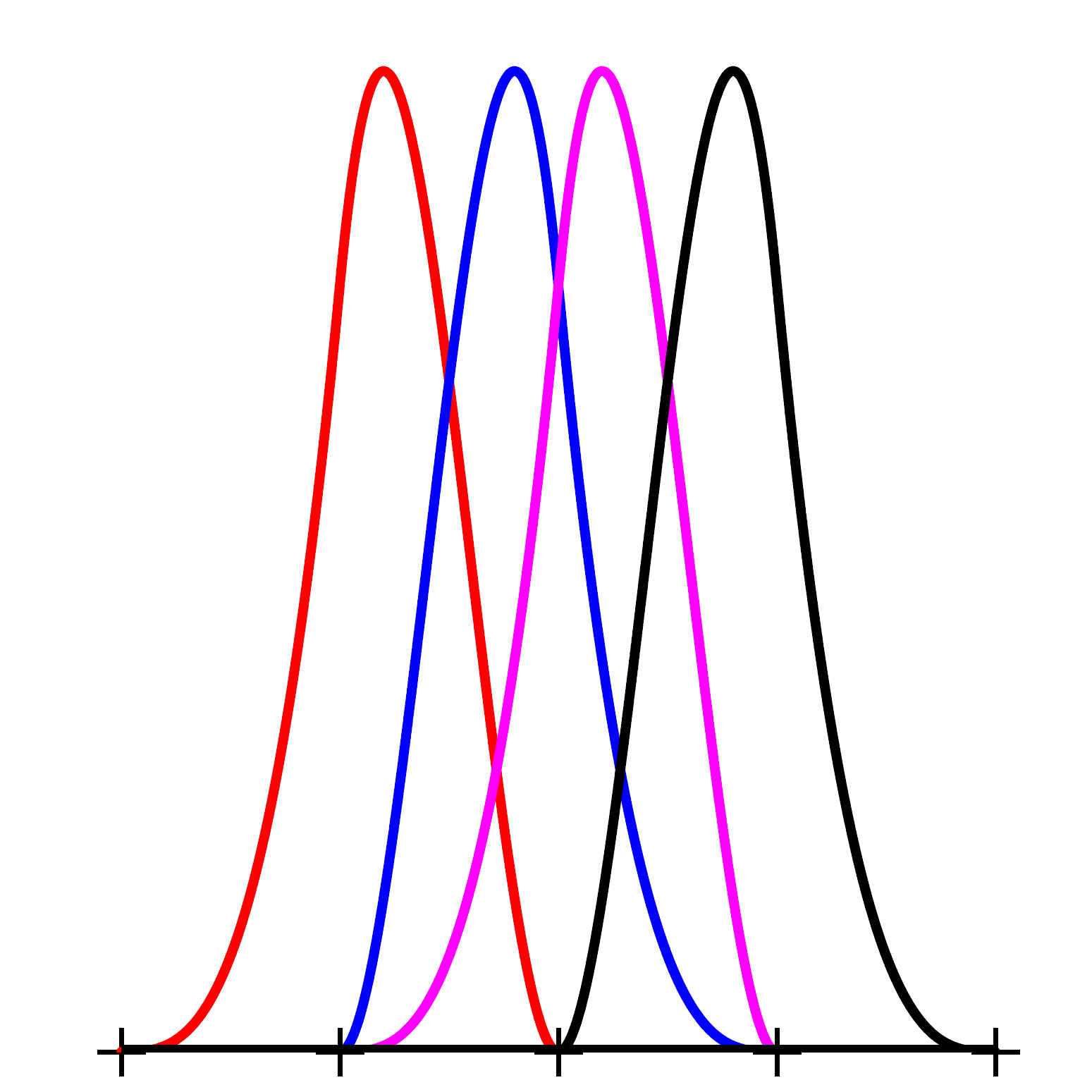}
 \put(0,70){\fcolorbox{gray}{white}{$\St_{3,1}^n$}}
    \put(28,-5){\small$\xt_i$}
    \put(21,90){\small$\Dt_j$}
    \put(2,10){\small$\XXt_n$}
	\end{overpic}
 \hfill
 \begin{overpic}[width=.32\columnwidth,angle=0]{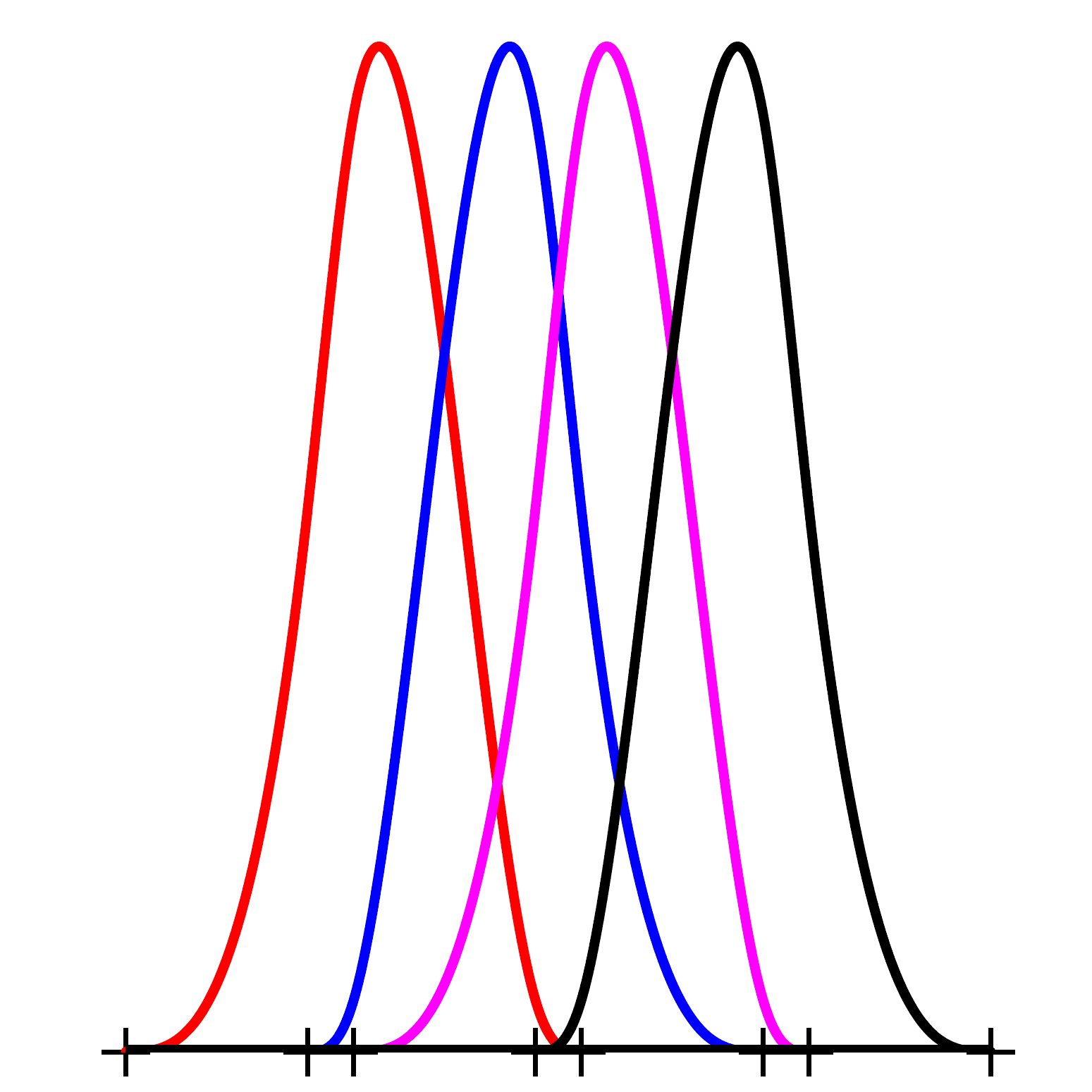}
 \put(-10,50){\huge$\rightarrow$}
	\end{overpic}
 \hfill
 \begin{overpic}[width=.32\columnwidth,angle=0]{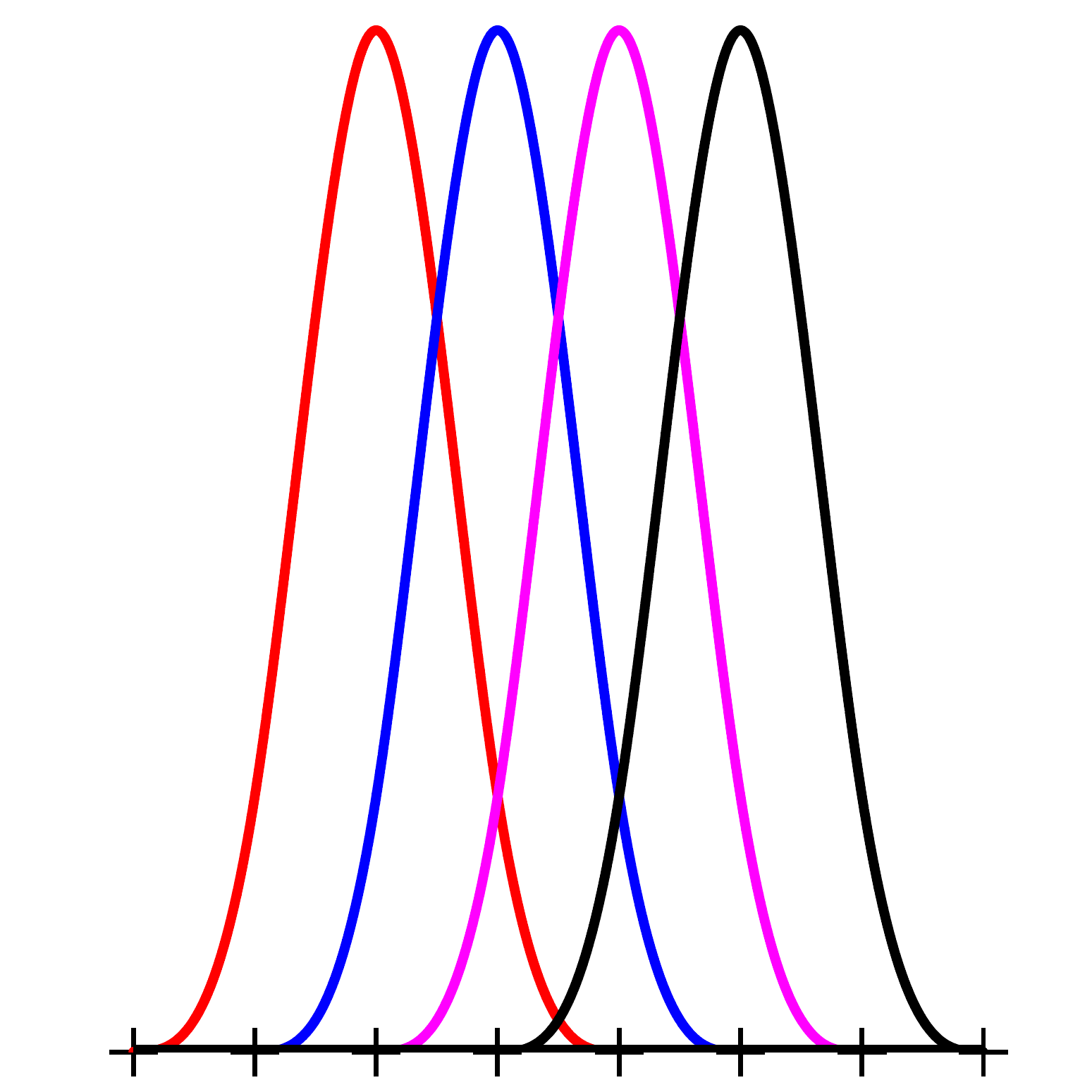}
 \put(80,70){\fcolorbox{gray}{white}{$S_{3,2}^N$}}
  \put(-10,50){\huge$\rightarrow$}
  \put(20,-5){\small$x_{2i-1}$}
  \put(20,90){\small$D_j$}
  \put(86,10){\small$\XX_N$}
	\end{overpic}\hfill \vrule width0pt\\
 \vspace{-5pt}
\Acaption{1ex}{Continuous transformation of cubic spline spaces. The uniform knot sequence with knots of multiplicity two
is transformed to a uniform knot sequence consisting of single knots. Four corresponding basis functions of
the source space $\St_{3,1}^n$, an intermediate space, and the target space $S_{3,2}^N$ are shown.}\label{fig:Trans}
 \end{figure}

\begin{rem}
The uniform knot sequences $\XXt_n$ and $\XX_N$ are, in this work, taken as the source and the target knot sequences, respectively.
Additionally, we consider the intermediate knot sequences generated as $\XXt_n$ continuously evolves to $\XX_N$, see Fig.~\ref{fig:Trans}.
Given a particular knot vector pattern, which may contain both single and double knots,
we also consider all the spaces of mixed continuities, not just (\ref{eq:familyC1}) and (\ref{eq:familyC2}).
\end{rem}

\begin{figure}[!tb]
\vrule width0pt\hfill
 \begin{overpic}[width=.82\textwidth,angle=0]{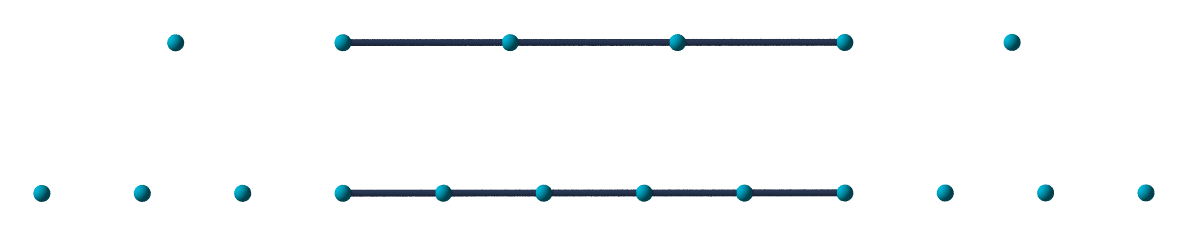}
    \put(-7,2){\fcolorbox{gray}{white}{$\XX_N$}}
    \put(-7,15){\fcolorbox{gray}{white}{$\XXt_n$}}
    \put(85.5,16){\vector(1,-1){10}}
    \put(85.5,16){\vector(1,-4){2.6}}
    \put(71.3,16){\vector(2,-3){7.2}}
    \put(71.3,16){\vector(0,-1){10.5}}
    \put(27,0){\small$a$}
    \put(71,0){\small$b$}
    \put(2,0){\small$x_{-3}$}
    \put(95,0){\small$x_{N+3}$}
    \put(25,12){\small$a=\xt_0$}
    \put(12,12){\small$\xt_{-1}$}
    \put(87,16){\small$\xt_{n+1}$}
    \put(29,17){$\overbrace{\mbox{\hspace{1.4cm}}}$}
    \put(29,3){$\underbrace{\mbox{\hspace{0.8cm}}}$}
    \put(35,20){\small$h$}
    \put(31,-2){\small$H$}
	\end{overpic}
 \hfill \vrule width0pt\\
 \vspace{-5pt}
\Acaption{1ex}{A target uniform knot vector with single knots $\XX_N$ for $N=5$ elements and an associated source knot vector $\XXt_n$
with all knots of multiplicity two are shown.
Both knot vector are of the same cardinality and have the same number of interior knots $i=4$.}\label{fig:Knots}
 \end{figure}

Consider now the $\St^{n}_{3,1}$ space.
Similarly to \cite{Nikolov-1996,Quadrature51-2014,Quadrature31-2014}, we work with the
non-normalized B-spline basis. To define the basis, we consider $\xt_0$ and $\xt_n$ as double knots and
extend our knot sequence $\XXt_n$ with two extra double knots outside the interval $[a,b]$ that we set to be
\begin{equation}\label{boundary_condition}
\xt_{-1} = \xt_{0} - h
\quad \textnormal{and} \quad
\xt_{n+1} = \xt_{n} + h,
\end{equation}
see Fig.~\ref{fig:Knots}.
The choice of $\xt_{-1}$ and $\xt_{n+1}$ yields particular integral values for the first and last
functions as discussed in (\ref{boundaryIntegral}).
Denote by $\DDt = \{\Dt_k\}_{k=1}^{2n+2}$ the basis of $\St^{n}_{3,1}$ where
\begin{equation*}
\renewcommand{\arraystretch}{1.2}
\begin{array}{rcl}
\Dt_{2k-1}(t) & = & [\xt_{k-2},\xt_{k-2},\xt_{k-1},\xt_{k-1},\xt_k](. - t)_{+}^{3}   \\
\Dt_{2k}(t)   & = & [\xt_{k-2},\xt_{k-1},\xt_{k-1},\xt_{k},\xt_k](. - t)_{+}^{3}, \\
\end{array}
\end{equation*}
where $[.]f$ stands for the divided difference and $u_{+} = \max(u,0)$ is the truncated
power function. 
Direct computation gives
\begin{equation}\label{interiorIntegral}
I[\Dt_k] = \frac{1}{4}\; \textnormal{for} \quad k = 3,4,\ldots,2n,
\end{equation}
where $I[f]$ stands for the integral of $f$ over the interval $[a,b]$.
We work with non-normalized basis functions, and
therefore the integrals above are independent on the knot sequence.
With the choice made in (\ref{boundary_condition}), direct integration gives
\begin{equation}\label{boundaryIntegral}
I[\Dt_1] =  I[\Dt_{2n+2}] = \frac{1}{16}
\quad \textnormal{and} \quad
I[\Dt_2] =  I[\Dt_{2n+1}]= \frac{3}{16}.
\end{equation}

Similarly to (\ref{boundary_condition}), we extend the knot sequence of $\XX_N$ by two triplets of single knots as
\begin{equation}\label{boundary_conditionC2}
x_{-k} = x_{0}- kH \quad \textnormal{and} \quad x_{N+k} = x_{N} + kH, \quad k=1,2,3.
\end{equation}
and define $\D = \{D_k\}_{k=1}^{2n+2}$ the basis of $S^{N}_{3,2}$ where
\begin{equation*}
\renewcommand{\arraystretch}{1.2}
\begin{array}{rcl}
D_{k}(t) & = & [\xt_{k-4},\xt_{k-3},\xt_{k-2},\xt_{k-1},\xt_{k}](. - t)_{+}^{3}  \quad k = 1, \dots, 2n+2. \\
\end{array}
\end{equation*}
We obtain 
\begin{equation}\label{eq:Areas}
I[\Dt_k] =  I[D_{k}], \quad k=4,\dots, 2n-1
\end{equation}
and the six boundary integrals (three only due to symmetry) are computed
directly by integrating $D_1$, $D_2$, and $D_3$ on $[a,b]$. These integrals
change during continuation and therefore have to be recomputed for various $t$.

\begin{figure}[!tb]
\vrule width0pt\hfill
 \begin{overpic}[width=.72\textwidth,angle=0]{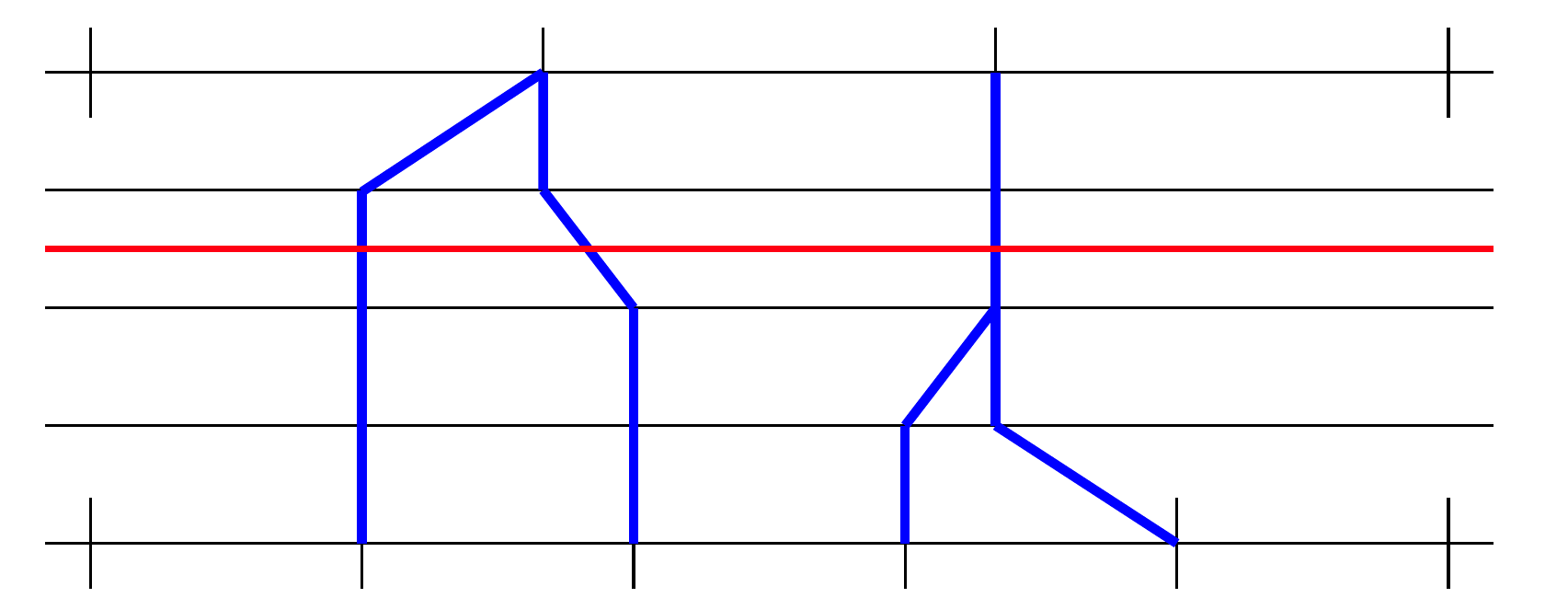}
 \put(-6,32){$\XXt_n$}
 \put(-6,3){$\XX_N$}
 \put(-6,20){$\XX(t)$}
 \put(30,35){$\xt_1$}
 \put(18,0){$x_1$}
 \put(41,0){$x_2$}
	\end{overpic}
 \hfill
 \begin{overpic}[width=.24\columnwidth,angle=0]{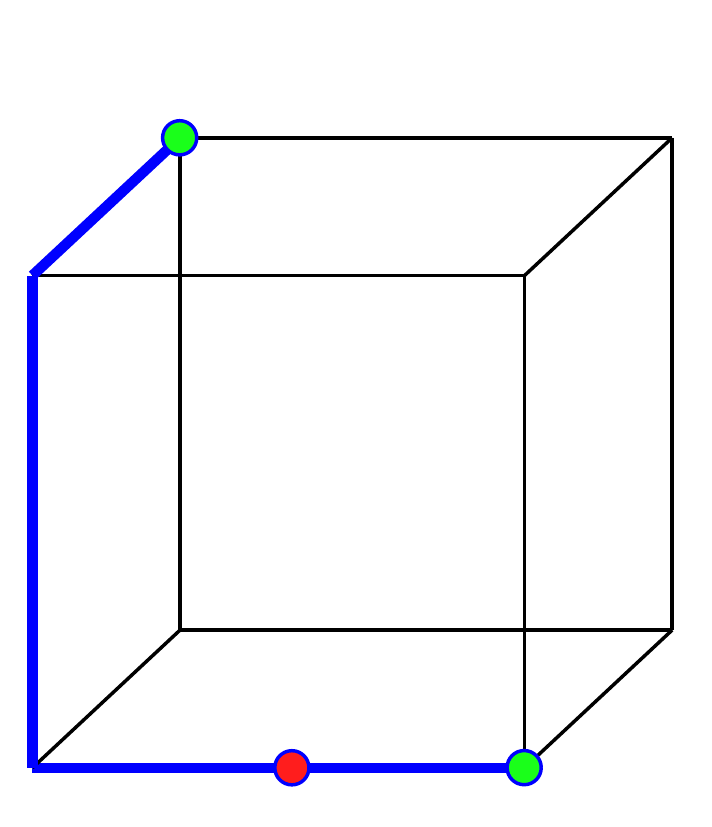}
 \put(90,60){$\R^{2n+4}$}
  \put(60,-5){\small$\XXt_n=\XX(0)$}
  \put(30,-5){\small$\XX(t)$}
  \put(10,90){\small$\XX_N=\XX(1)$}
	\end{overpic}\hfill \vrule width0pt\\[3ex]
\vrule width0pt\hfill
 \begin{overpic}[width=.72\textwidth,angle=0]{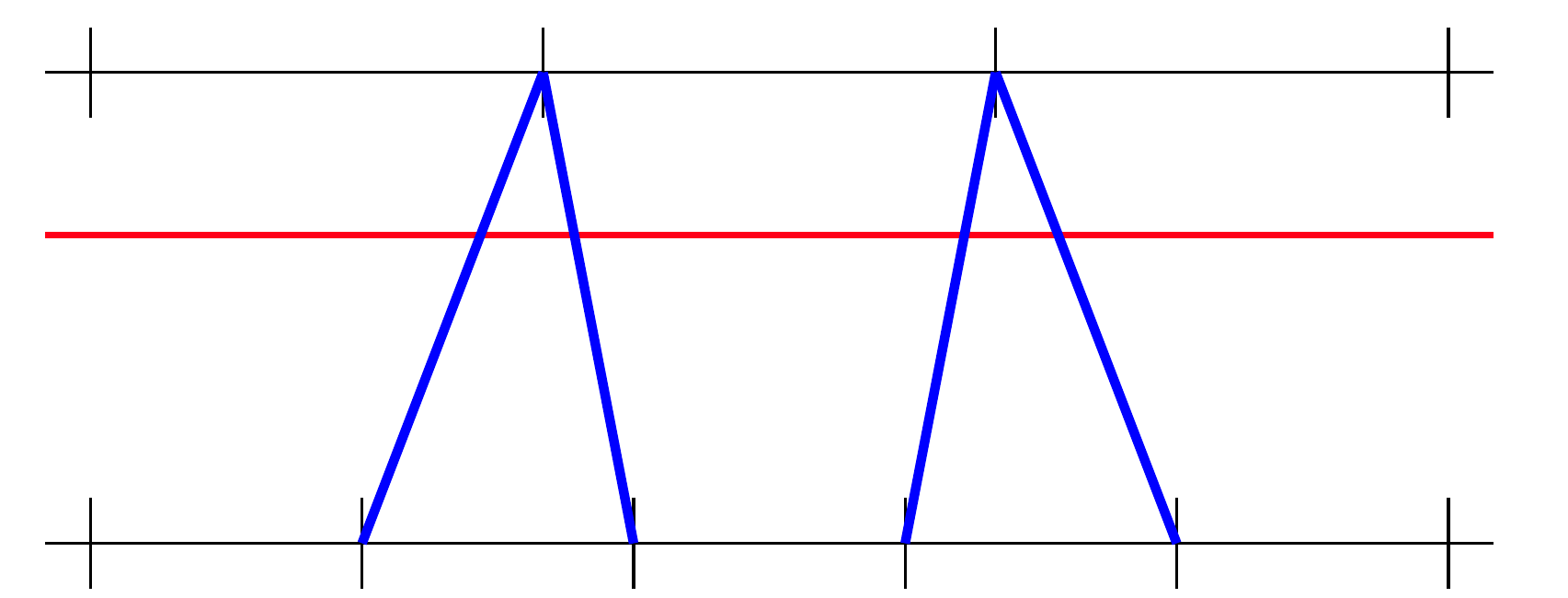}
\put(-6,32){$\XXt_n$}
 \put(-6,3){$\XX_N$}
 \put(-6,20){$\XX(t)$}
 \put(30,35){$\xt_1$}
 \put(18,0){$x_1$}
 \put(41,0){$x_2$}
	\end{overpic}
 \hfill
 \begin{overpic}[width=.24\columnwidth,angle=0]{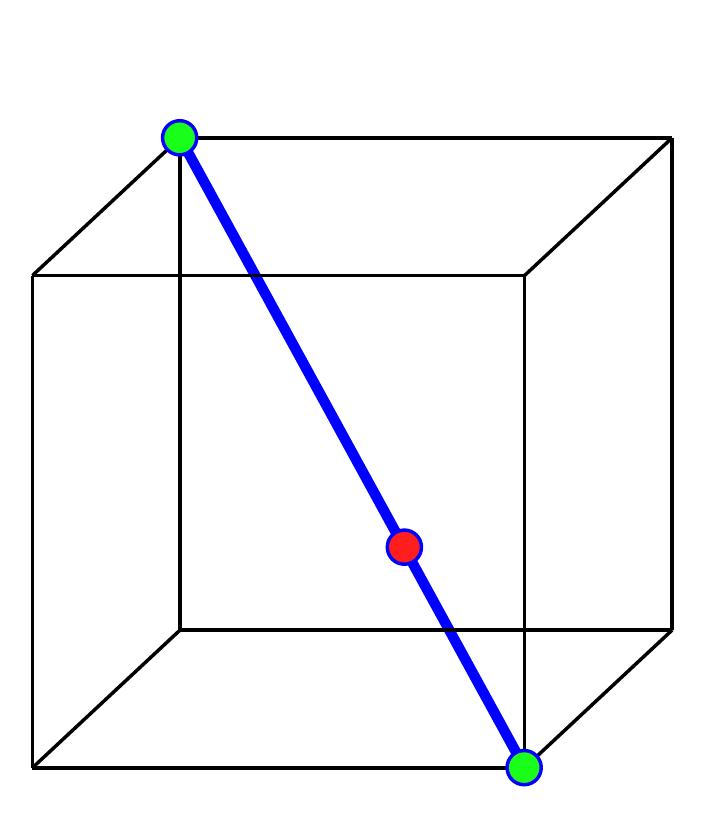}
 \put(22,30){$\XX(t)$}
\end{overpic}\hfill \vrule width0pt\\[-2ex]
\Acaption{1ex}{Continuous transformation of knots sequences.
The source knot sequence $\XXt_n$ with double knots is changed over time and is transformed into the target uniform
knot sequence with single knots, $\XX_N$.
The transformation is metaphorized as a path in $2n+4$ dimensional space of free knots.
Two particular paths are shown: whereas one transformation changes always only a single knot (top),
the geodesic transformation (bottom) generates the shortest path in the space of knots
and consequently the shortest transformation from $\XXt_n$ into $\XX_N$.}\label{fig:Transknots}
 \end{figure}

Now we consider a continuous transition between $\St_{3,1}^n$ and $S_{3,2}^N$, see Fig.~\ref{fig:Trans}.
Since the transition of the spline spaces
is governed by the transformation of the corresponding knot vectors,
consider the following mapping
\begin{equation}
\XXt_n \rightarrow \XX_N,
\end{equation}
including the six outer knots defined in (\ref{boundary_condition}) and (\ref{boundary_conditionC2}),
see Fig.~\ref{fig:Knots}.
Due to the fact that we work with non-normalized basis functions,
Eqns. (\ref{eq:Areas}) remain unchanged.
The total number of knots is $2n+6$, but since the two boundary knots are constrained to stay fixed,
there is $2n+4$ free knots.
%
The transformation can be conceptualized
as a curve between $\XXt_n$ and $\XX_N$, two points in $\R^{2n+4}$, see Fig.~\ref{fig:Transknots}.
There exist infinitely many paths connecting the source and target knot vectors. In particular, we analyze two specific
knot transformations:
the first one sequentially changes only one knot whilst all the
others remain unchanged; the second one simultaneously spreads all free knots from the source position
to the target one in a linear fashion. Such a transition can be seen as a diagonal straight line connecting $\XXt_n$ and $\XX_N$,
i.e., a geodesic path when under the Euclidean metric on the vector of free knots,
while the first one corresponds to a path along the edges of a $(2n+4)$-dimensional hypercube.

\begin{rem}
Alternatively, one can extend $\XXt_n$ and $\XX_N$ by adding different six knots than those in
(\ref{boundary_condition}) and (\ref{boundary_conditionC2}). For example, one can extend them in the open knot fashion
by setting the multiplicity four to both boundary knots.
In that case, the six (three due to symmetry) initial boundary integrals have to be computed accordingly.
\end{rem}

\section{Gaussian quadrature rules via homotopy continuation}\label{sec:HomotopyQuad}

In this section, we use polynomial homotopy continuation (PHC) to generate
Gaussian quadrature rules for families of cubic spline spaces above arbitrary knot sequences.

The key observation is the following.
Assume you have an exact and optimal quadrature rule for a spline space $\St$ 
defined above a knot vector $\XXt$ and consider a continuous transformation of $\XXt$ as a function of time, i.e., $\XXt(t)$.
As a quadrature rule is a solution of a particular
polynomial system defined above $\XXt$, the quadrature
nodes and weights also change continuously. We use
the PHC framework, to take an exact and optimal rule as the initial
root of a certain polynomial system and trace this high-dimensional point
(``follow the path of the root'' in the PHC terminology, see \cite{Wampler-2005})
while changing continuously $\XXt$. This idea is in accordance with the result of Micchelli and Pinkus \cite{Micchelli-1977}
which states that there exists an optimal quadrature rule with $m$ nodes such that
\begin{equation}\label{eq:Micchelli}
     d + i + 1 = 2 m
\end{equation}
where $d$ is the spline degree and $i$ is the number of interior knots, counted including their multiplicities.
Therefore $m$ stays fixed as long as the number of interior knots remains constant.
Homotopy continuation transfers the optimal quadrature rule from one space to another
because the number of interior knots remains unchanged.


In this paper, we want to derive optimal quadrature rules for spaces of $C^2$ cubic splines,
so for uniform knot vectors we have $S_{3,2}^N$ as our target spaces, see (\ref{eq:familyC2}).
The multiplicity of the interior knots is one and therefore Eq.~(\ref{eq:Micchelli}) requires $N$ to be odd.
Further, we have $n=[\frac{N}{2}]+1$ and set $\St_{3,1}^n$ (\ref{eq:familyC1}) as our source spaces.
For these spaces, unique quadrature rules that are explicit were derived in \cite{Nikolov-1996}.
The framework we present is general and one can derive new rules by applying PHC to any other
appropriate source space where an optimal quadrature rule is known,
e.g., $C^1$ quintic splines with uniform knot sequences \cite{Quadrature51-2014}.

\subsection{Gaussian quadrature formula}\label{ssec:Homo}

\begin{figure}[!tb]
\vrule width0pt\hfill
 \begin{overpic}[width=.45\textwidth,angle=0]{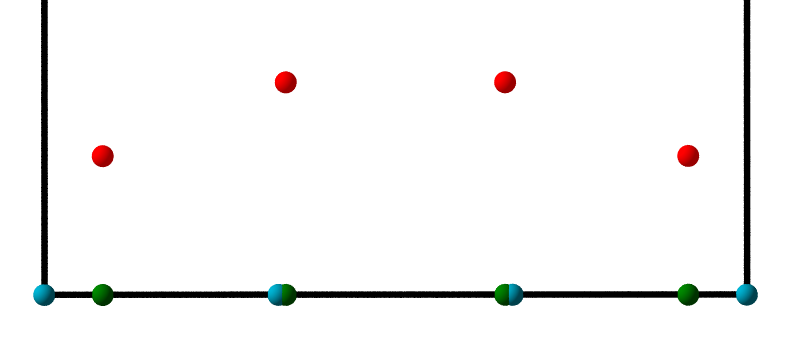}
 \put(20,50){\fcolorbox{gray}{white}{source space $\St_{3,1}^n$}}
 \put(25,17){\fcolorbox{gray}{white}{source rule $\Qt$}}
    \put(2,-1){\small$a$}
    \put(11,-1){\small$\taut_1$}
    \put(39,31){\small$[\taut_2,\omegat_2]$}
    \put(30,-1){\small$\xt_1$}
    \put(65,-1){\small$\xt_2$}
    \put(37,-1){\small$\taut_2$}
    \put(95,-1){\small$b$}
	\end{overpic}
 \hfill
 \begin{overpic}[width=.45\textwidth,angle=0]{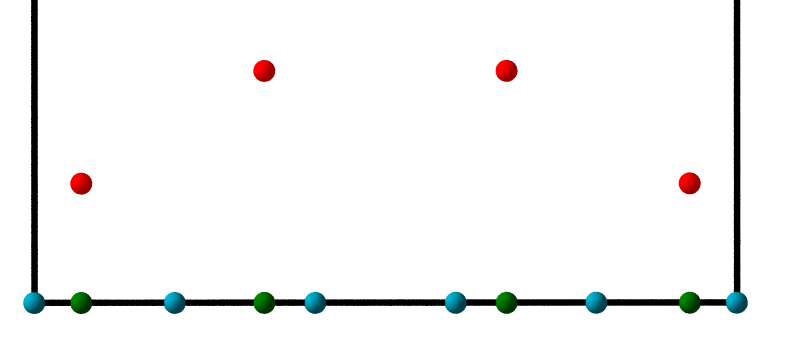}
    \put(20,50){\fcolorbox{gray}{white}{target space $S_{3,2}^N$}}
    \put(25,17){\fcolorbox{gray}{white}{target rule $\Q$}}
    \put(-10,50){\huge$\rightarrow$}
    \put(2,-1){\small$a$}
    \put(20,-1){\small$x_1$}
    \put(75,-1){\small$x_4$}
    \put(32,-1){\small$\tau_2$}
    \put(37,35){\small$[\tau_2,\omega_2]$}
    \put(90,-1){\small$b=x_5$}
	\end{overpic}\hfill \vrule width0pt\\
 \vspace{-5pt}
\Acaption{1ex}{For a desired (target) space $S_{3,2}^N$ with $N=5$ elements on $[a,b]$, an associated source space $\St_{3,1}^n$, $n=3$, is built (left).
The source knot sequence (blue) is uniform with two interior knots $\xt_1$ and $\xt_2$, each of multiplicity two.
The target knot sequence has four single knots, so the total number of interior knots is unchanged; $i=4$ in (\ref{eq:Micchelli}).
Therefore both optimal rules require $m=n+1$ nodes (green) and weights (red).
The target rule (right) is derived via homotopy continuation, see Algorithm~\ref{algor}.}\label{fig:SourceTarget}
 \end{figure}

We set our source space as $\St^{n}_{3,1}$, see (\ref{eq:familyC1}) and
know, according to (\ref{eq:Micchelli}), that $m=n+1$. The optimal
source quadrature rule is of the form
\begin{equation}\label{quadratureS}
\Qt_a^b[f] = \int_{a}^{b} f(t) \mathrm{d}t = \sum_{i=1}^{n+1} \omegat_i f(\taut_i)
\end{equation}
and the nodes and weights can be computed by a recursion derived by Nikolov \cite{Nikolov-1996}, see Fig.~\ref{fig:SourceTarget}.
Due to the equal dimensions of $\St^{n}_{3,1}$ and $S^{N}_{3,2}$, the target rule requires the same number of nodes
and therefore
\begin{equation}\label{quadrature}
\Q_a^b[f] = \int_{a}^{b} f(t) \mathrm{d}t = \sum_{i=1}^{n+1} \omega_i f(\tau_i).
\end{equation}

During the continuation, we transform the spline space $\St^{n}_{3,1}$ to $S^{N}_{3,2}$
and accordingly the optimal rule $\Qt \rightarrow \Q$.
Therefore $\Q$, represented by the its nodes and weights, is a function of $t$.
To simplify notation, if no ambiguity is imminent,
we omit the time parameter and write $\tau_i$ instead of $\tau_i(t)$.
The source rule is $\Qt=\Q(0)$ and the target rule is $\Q=\Q(1)$.
We denote by $\Q$ the rule (linear operator) and later in Section~\ref{ssec:IrregNP}
use the symbol $\br$ for a $2m$-dimensional point,
but they both refer to the same set of optimal nodes and weights.
Before we proceed to our homotopy continuation setting, we need to establish notation
that unifies the source and the target quadrature rules.

Let $\XX_N = (a=x_0,x_1,\dots,x_{N-1},x_{N} = b)$ be a knot vector consisting of $N$ subintervals,
some of the subintervals may degenerate to zero length
in the case of a double knot and let $\{\tau_1,\dots, \tau_m\}$ be the quadrature points.
We define a \emph{nodal pattern} $\bp$ of the quadrature rule $\Q$
on $\XX_N$ as an $N$-dimensional vector where the $i$-th coordinate specifies the number of quadrature nodes inside
the $i$-th sub-interval. We say $\bp$ is \emph{regular} if no node coincides with a knot.

%

\begin{ex}
Consider the Gaussian quadrature rule of Nikolov \cite{Nikolov-1996} for (\ref{eq:familyC1})
in the case when $n$ is odd, $n>1$. The rule says every sub-interval contains exactly one node except the middle one
which contains two, see Fig~\ref{fig:SourceTarget} for $n=3$.
Considering the multiplicities, the $\St_{3,1}^n$ contains $n$ sub-intervals of non-zero length
and $(n-1)$ degenerated subintervals (double knots). The nodal pattern is then
\begin{equation}\label{eq:NP31odd}
\begin{array}{cccc}
\bp =  & (\underbrace{1,0, \dots,1,0,} & \underbrace{2,} & \underbrace{0,1,\dots,0,1} )\\
& n-1 & 1 & n-1
\end{array}
\end{equation}
and the number of nodes is exactly $m=n+1$.
In the case when $n$ is even, the middle node coincides with the middle knot. In such a case, we
say the nodal pattern is {\em irregular}. We will discuss these patterns later in Section~\ref{ssec:IrregNP}.
\end{ex}

The nodal patterns describe the layout of the quadrature nodes with respect to the knots
and, as long as the pattern does not change, the pattern of the polynomial system
to solve remains unchanged. The crucial part of any optimal quadrature rule is to derive
a correct nodal pattern. The result of Micchelli and Pinkus \cite{Micchelli-1977} states
the number of optimal nodes and also the knot span in which each particular
node lies. However, one cannot conclude the nodal pattern for $S_{3,2}^N$ from their result.
In our approach once we know the nodal pattern of the source rule,
then we know the initial polynomial system and its root.

\subsection{Homotopic setting}\label{ssec:Setting}

We are now ready to apply the PHC framework to the quadrature problem.
The vector of unknowns consists of the quadrature nodes and weights
$$
\bx = (\tau_1,\dots, \tau_m, \omega_1,\dots, \omega_m), \quad \bx \in \R^{2m},
$$
our source polynomial system $\FFt$ expresses that the source rule $\Qt$ exactly
integrates the source basis $\DDt$, that is,
\begin{equation}\label{eq:IniSystem}
\Qt_a^b(\Dt_i) = I[\Dt_i], \quad i = 1,\dots, 2m
\end{equation}
and the source root $\br$ that solves (\ref{eq:IniSystem}) is the quadrature rule of Nikolov \cite{Nikolov-1996}.
The domain $\Omt \subset \R^{2m}$ is generated as follows. For this quadrature rule we know
the nodal pattern, e.g., for $n$ odd, $n=2k+1$, we have
(\ref{eq:NP31odd}), therefore
\begin{equation}\label{eq:OmegaNodes}
(\tau_1, \dots, \tau_m) \in [\xt_0,\xt_1]\times \dots \times [\xt_k, \xt_{k+1}] \times [\xt_k, \xt_{k+1}] \times \dots \times [\xt_{n-1}, \xt_n]
\end{equation}
and for the weights we use (a very rough) range $[0,b-a]$. Combined together, the source domain is
\begin{equation}\label{eq:OmegaT}
\begin{array}{cccc}
\Omt = & \underbrace{[\xt_0,\xt_1] \times \dots \times [\xt_{n-1}, \xt_n]} & \times & \underbrace{[0,b-a] \times \dots \times [0,b-a]}.\\
& m & & m
\end{array}
\end{equation}

\begin{rem}
Eq.~(\ref{eq:OmegaT}) gives a very loose bound on $\omega_i$, $i=1,\dots,m$.
However, this part of $\Omt$ is not affected by the nodal pattern and thus less important
since it does not influence the change of the pattern of the polynomial system. A tighter bound could
eventually serve as a better stopping criterion in cases, where the optimization is required to
keep the root inside the domain (see Section~\ref{ssec:Impl}) but 
we expect this would bring only marginal computational gains.
\end{rem}

There are several differences between the homotopy framework we propose and the classical
setting \cite{Wampler-2005}:
\begin{itemize}
\item There is only one root that we follow: the quadrature rule $\Qt$ of $\St_{3,1}^n$. 
Moreover, we know the quadrature rules (the source rule, the intermediate ones and the target rule) are {\em unique}
for geodesic knot transformations.
This fact is a great advantage as there is no danger of numerical instabilities like jumping from one traced path to another
as is the case in classical homotopy continuation methods.
\item The systems do not require linear blending as in \cite{Wampler-2005}. In our case, the continuous
transition between the systems is governed by the change of the knot vector.
\item The domain $\Omt$ also changes in time because the range of every node
is determined by the knot positions and these vary.
%
\item The system is only {\em locally} polynomial. When a node leaves its subinterval,
the nodal pattern of the quadrature rule changes,
a new polynomial system has to be built and solved.
\end{itemize}

The last issue is important. The most difficult part of any optimal spline quadrature
rule is to derive the correct layout of the nodes (the nodal pattern).
Since the splines are piece-wise polynomials, the systems that we build and solve are
also only ``piece-wise polynomial'', i.e. locally polynomial and
the main difficulty is to select the ``right pieces''. However, homotopy continuation leads us
to the correct nodal pattern automatically, by tracing down the trajectory of
optimal quadrature rules as the knot pattern evolves. 
That is, the generated quadrature rules are exact (up to machine precision, see
Section~\ref{sec:Ex}) and optimal for all the intermediate steps of the transformation.

\subsection{Irregular nodal patterns}\label{ssec:IrregNP}

\begin{figure}[!tb]
\vrule width0pt\hfill
 \begin{overpic}[width=.5\textwidth,angle=0]{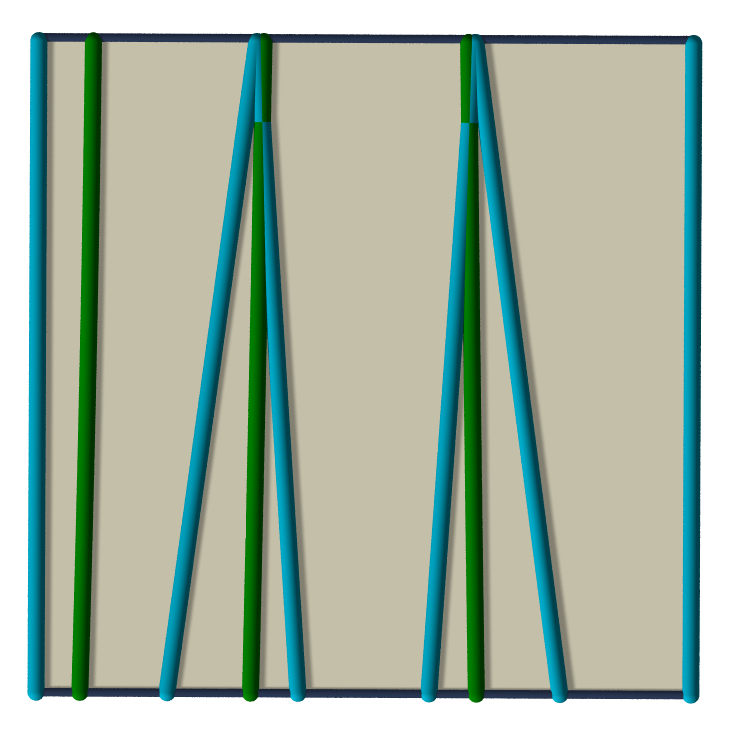}
    \put(2,-1){\small$a$}
    \put(18,-1){\small$x_1$}
    \put(39,-1){\small$x_2$}
    \put(55,-1){\small$x_3$}
    \put(32,-1){\small$\tau_2$}
    \put(62,-1){\small$\tau_3$}
    \put(95,-1){\small$b$}
    \put(-45,87){$\bp = (1,0,2,0,1)$}
    \put(-45,42){$\bp = (1,1,0,1,1)$}
      \thicklines
    \put(-3,82.5){\line(1,0){108}}
    \put(105,95){\vector(0,-1){91}}
    \put(107,95){\small$t=0$}
    \put(107,80){\small$t=t_1$}
    \put(107,2){\small$t=1$}
    \put(-3.5,5){\rotatebox{90}{$\overbrace{\color{white}{aaaaaaaaaaaaaaaaaaaaaaaaa}}$}}
    \put(-3.5,83){\rotatebox{90}{$\overbrace{\color{white}{aaaa}}$}}
    \put(87,87.5){\fcolorbox{gray}{white}{\small$I[D_2]=\omega_1 D_2(\tau_1)$}}
    \put(87,42){\fcolorbox{gray}{white}{\small
        $
        \begin{array}{ccc}
        I[D_2] & =  & \omega_1 D_2(\tau_1)\\
               & +  & \omega_2 D_2(\tau_2)
        \end{array}
        $
        }}
	\end{overpic}
 \hfill \vrule width0pt\\[2ex]
 \vrule width0pt\hfill
 \begin{overpic}[width=.99\textwidth,angle=0]{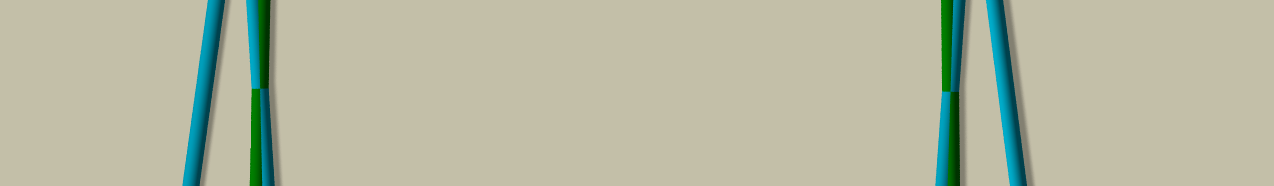}
 \thicklines
 \put(0,7.5){\line(1,0){100}}
	\end{overpic}
 \hfill \vrule width0pt\\
 \vspace{-3ex}
\Acaption{1ex}{Nodal pattern change.
Top: evolution of the nodes (green) and knots (blue) is shown. At time instant $t_1$, the source nodal pattern changes which corresponds
to the fact that $\tau_2$ leaves $(x_2, x_3)$ by crossing the knot $x_2$. Due to symmetry, $\tau_3$ leaves $(x_2, x_3)$ and
the new nodal pattern is $\bp = (1,1,0,1,1)$. The new nodal pattern requires a change of the system (\ref{eq:IniSystem}),
e.g., $\tau_2$ starts affecting the basis function $D_2$ and the second equation of the system must be updated (right).
Bottom: A zoom-in on the critical instant $t=t_1$.}\label{fig:NodalPattern}
 \end{figure}

We now discuss the changing the nodal patterns, see Fig.~\ref{fig:NodalPattern}.
This is the situation when
some (one or more) node(s) crosses the boundary of its subinterval.
That is, the node overlaps with a knot and the nodal pattern of the quadrature becomes irregular.
We recall that $\Om(t)$ is a hypercube in $\R^{2m}$ and the system (\ref{eq:IniSystem})
can be seen as an intersection problem of $2m$ hypersurfaces inside $\Om(t)$,
which has at time $t=0$ only one root $\br(0)$.
During the continuation, as time evolves, $\Om(t)$ is changing and so does the root $\br(t)$
that is being traced.

The situation when a node coincides with a knot corresponds to $\br(t)$
leaving $\Om(t)$, i.e., there exist some $t_1 \in [0,1]$ such that $\br(t_1)$
is on the boundary of the hypercube $\Om(t_1)$, $\br(t_1)\in \overline{\Om}(t_1)$.
For simplicity, we assume that $\br(t_1)$ lies on a hyperface of $\Om(t_1)$ which
corresponds to the fact that only {\em one} node becomes a knot, $\tau_i(t_1) = x_j(t_1)$.
Let the current sub-interval where $\tau_i(t)$ lies in for $t<t_1$ be $(x_{j-1}(t),x_{j}(t))$.
We know the hyper-face that is being crossed, $\tau_i(t_1) = x_j(t_1)$, so the algorithm changes the nodal pattern
by switching the sub-interval of $\tau_i$ from $(x_{j-1},x_{j})$ to $(x_{j},x_{j+1})$.
Then the system (\ref{eq:IniSystem}) is updated accordingly.

If $\br(t_1)$ lies on a hyper-edge, i.e., two (or more) nodes become knots at the same time,
the situation is similar because for every node we know the current interval that the node leaves
and the new one that it enters.


\subsection{Quadrature rule tracing}\label{ssec:Tracing}

The goal is to trace down a curve $\br(t)\in \R^{2m}$, $t\in[0,1]$ knowing the initial
point $\br(0)$, for example the quadrature rule of Nikolov \cite{Nikolov-1996}.
We know the source and the target knot sequences, $\XXt_n$ (\ref{eq:XXt}) and $\XX_N$ (\ref{eq:XX}), respectively.
As the number of non-degenerate intervals varies from $n$ to $N$, we simplify the notation by omitting the subscripts
and write $\XX(0):=\XXt_n$ and $\XX(1)=\XX_N$. We recall that the number of interior knots and total knots
remains constant. We further discretize time and build $M-1$
intermediate knot vectors
\begin{equation}\label{eq:InterKnots}
\{\XX(t_i) \}_{i=0}^M \quad (0=t_0,t_1, \dots, t_M = 1)
\end{equation}
and consequently discretize $\br(t)$ as a polyline $\{\br^i \}_{i=0}^M$.
We follow two different types of generation of the intermediate knot sequences, see Section~\ref{sec:Splines}.
The geodesic path, where each knot moves linearly in $t$ to its target position,
and along-the-edge paths, where only a single knot moves at each time instant, see Fig.~\ref{fig:Transknots}.
The first path is unique; the number of the paths of the second kind equals the number of ``bottom-to-top'' paths
on a hypercube in $\R^{2n+4}$, $(2n+4)!$.
The latter paths generate spline spaces with mixed continuities (various knot multiplicities).

\begin{rem}
The source quadrature rule is unique and so is the target rule.
Therefore the target rule should get derived independently on the path between $\XX(0)$ and $\XX(1)$.
However, to the best of our knowledge, there are no theoretical results on the quadrature rules for splines with {\em mixed} continuities.
The results of Micchelli and Pinkus \cite{Micchelli-1977} apply only for families of splines with uniform continuity.
Our target quadrature rule gets derived independently of the path chosen, by parsing these spaces with mixed continuities
as seen later in Section~\ref{sec:Ex}. This result provides numerical evidence for the existence of optimal rules 
for these mixed continuity spaces.
\end{rem}

\subsection{Implementation}\label{ssec:Impl}

\begin{algorithm}[!t]\normalsize \caption{\hfill { {\tt
GaussianQuadrature}$([a,b],N)$} }\label{algor}
\begin{algorithmic}[1]
\medskip
\STATE \textbf{INPUT}: compact interval $[a,b]$ and odd
number of uniform segments $N$\\
\vspace{-0.35cm}
 \begin{tabular}{c}
  \makebox[11cm]{}\\ \hline
 \end{tabular}
\vspace{0.15cm}
\STATE $n=[\frac{N}{2}]+1$;
\STATE $\br^0:=$ the source rule $\Qt$ on $[a,b]$; 
\STATE build $\FF^0(\bx)=\mathbf 0$ from (\ref{eq:IniSystem}) over $\Om^0$ (\ref{eq:OmegaT});
\STATE build $\{\XX(t_i)\}_{i=0}^{M}$, Section~\ref{sec:Splines};
\FOR{$i=1$ to $M$} 
\STATE build $\FF^{i}(\bx)=\mathbf 0$ over $\Om^{i}$;
\STATE $\br^{i}_{opt}:=$ solution of $\FF^{i}(\bx)=\mathbf 0$ by multivariate Newton-Raphson \\ with initial guess $\br^{i-1}$;
\IF{$\br^{i}_{opt} \in \Om^i$}
    \IF{$\|\br^{i}_{opt}\| < \varepsilon$}
    \STATE $\br^i:=\br^{i}_{opt};$ \hfill {\tt /* the quadrature rule for $\XX(t_i)$ */}
    \ELSE
    \STATE subdivide $\XX(t_{i-1})$ and $\XX(t_{i})$; \hfill {\tt /* finer stepsize */}
    \STATE $M:=M+1;$ and go to line 7;
    \ENDIF
\ELSE
\STATE update $\FF^i$ and $\Om^i$ and go to line 7; \hfill {\tt /* nodal pattern changes */}
\ENDIF
\ENDFOR
\STATE \textbf{OUTPUT}: $\br^M$, the Gaussian quadrature rule for $S_{3,2}^N$ on $[a,b]$;
\end{algorithmic} \end{algorithm}

In our implementation, the time-discretization of the intermediate knot vectors is uniform.
In the case of the geodesic knot transformation, we
sample uniformly the diagonal (source-target) of the hypercube, see Fig.~\ref{fig:Transknots}.
In the case of the along-the-edge transformation, we sample uniformly every hyperedge.

However, one can apply a non-uniform stepsize, e.g. this stepsize could get increased in cases where the nodes
differ from the knots by more than a certain threshold. Reversely, a finer modification
of the knot vector would be beneficial in cases when a node approaches a knot, i.e., in cases that
precede a change of the nodal pattern. Such an analysis
goes beyond the scope of this paper.
We set $M=200$ for the geodesic path and $M=20(N-1)$ for along-the-edge paths.

The tracing algorithm sequentially transforms the knot sequence $\XX(t_i)$ and updates the system~(\ref{eq:IniSystem}),
$\FF^{i}(\bx)=\mathbf 0$ considered over the domain $\Om^{i}$.
The exact quadrature rule from the previous iteration $\br^{i-1}\in \R^{2m}$ is taken as an initial guess and a multivariate Newton-Raphson iterative scheme
is applied to get the root $\br^i_{opt}$. If $\br^i_{opt}$ lies in the current domain $\Om^i$, the next knot sequence $\XX(t_{i+1})$
is taken with $\br^{i}:=\br^i_{opt}$ and the algorithm proceeds iteratively.
When $\br^i_{opt} \notin \Om^{i} $, this implies that the nodal pattern has changed as explained in Section~\ref{ssec:IrregNP}.
In these instances, the algorithm swaps to a different polynomial system defined above $\Om^{new}$.
Then $\br^{i}:=\br^i_{opt}$ and the Newton-Raphson optimization is applied again
with the new system $\FF^{new}(\bx)=\mathbf 0$.
If the optimization fails, i.e., if $\|\br^i_{opt}\|>\varepsilon$ an additional intermediate knot sequence is inserted between
$\XX(t_{i-1})$ and $\XX(t_{i})$ and the original point $\br^{i-1}$ and the original system $\FF^{i}(\bx)=\mathbf 0$ are considered again.
Otherwise the nodal pattern is updated and $\Om^{i}:=\Om^{new}$. The numerical threshold was set $\varepsilon = 10^{-16}$.
The whole procedure is summarized in Algorithm~\ref{algor}.

\begin{rem}
Algorithm~\ref{algor} requires the number of elements $N$ to be odd. This is because the associated source space
has all interior knots of multiplicity two and therefore $N=n+n-1$. This is in accordance with the result of Micchelli (\ref{eq:Micchelli})
that states the optimal rule for $S^N_{3,2}$ is guaranteed only for odd $N$.
\end{rem}

\section{Numerical examples}\label{sec:Ex}

\begin{figure}[!tb]
\vrule width0pt\hfill
 \begin{overpic}[width=.89\textwidth]{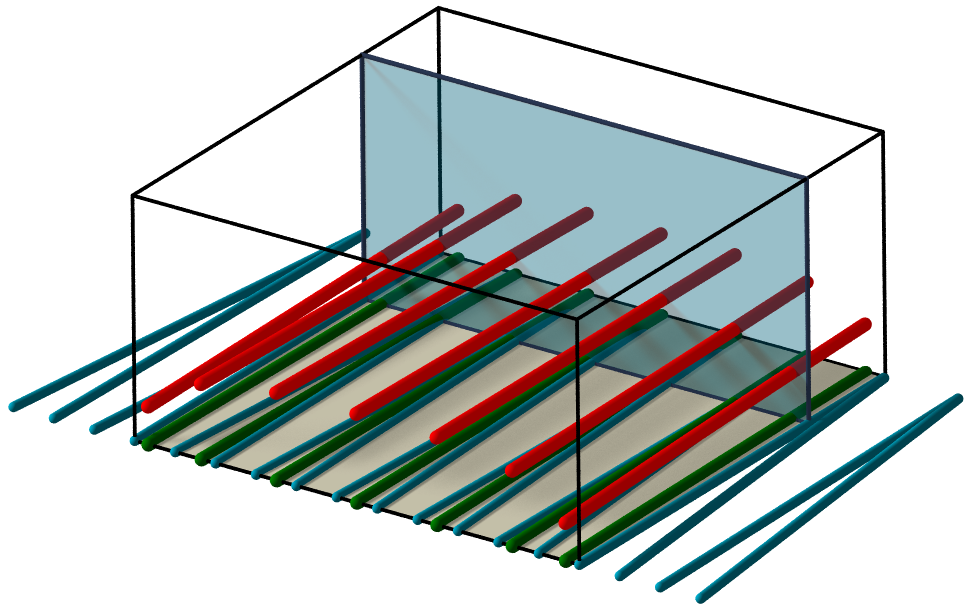}
 \put(-5,60){(a)}
    \put(11,15){\small$a$}
    \put(57,2){\small$b$}
    \put(6,44){\small$t=1$}
    \put(26,58){\small$t=0.25$}
    \thicklines
    \put(105,20){\vector(-4,-3){20}}
    \put(96,10){\small$t$}
    \put(105,20){\vector(0,1){25}}
    \put(107,35){\small$\omega$}
    \put(92,50){\small$[0,b,0.5]$}
	\end{overpic}
 \hfill \vrule width0pt\\[-1.ex]
 \vrule width0pt\hfill
 \begin{overpic}[width=.34\columnwidth,angle=0]{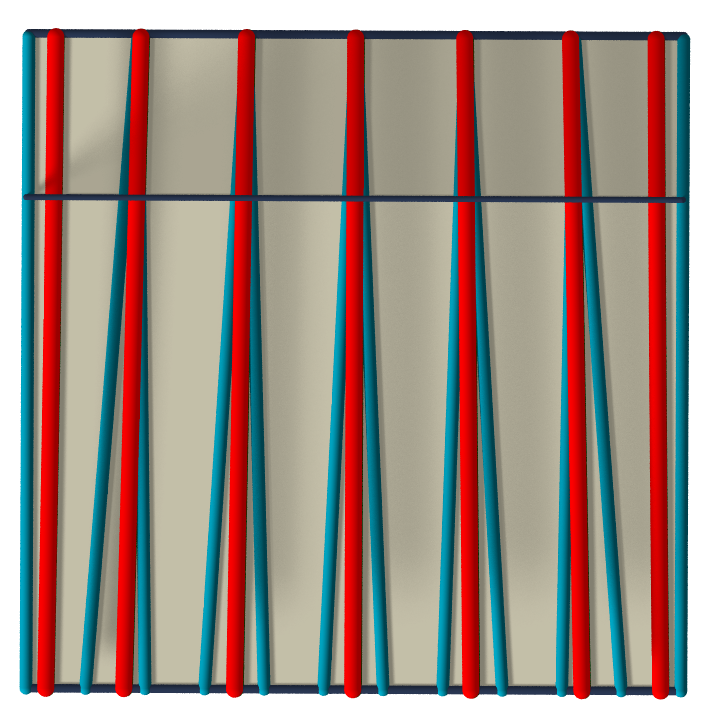}
 \put(0,100){(b)}
 \put(103,100){(c)}
\put(30,-5){\small$t=1$}
\put(30,100){\small$t=0$}
\end{overpic}
\hfill
 \begin{overpic}[width=.62\columnwidth,angle=0]{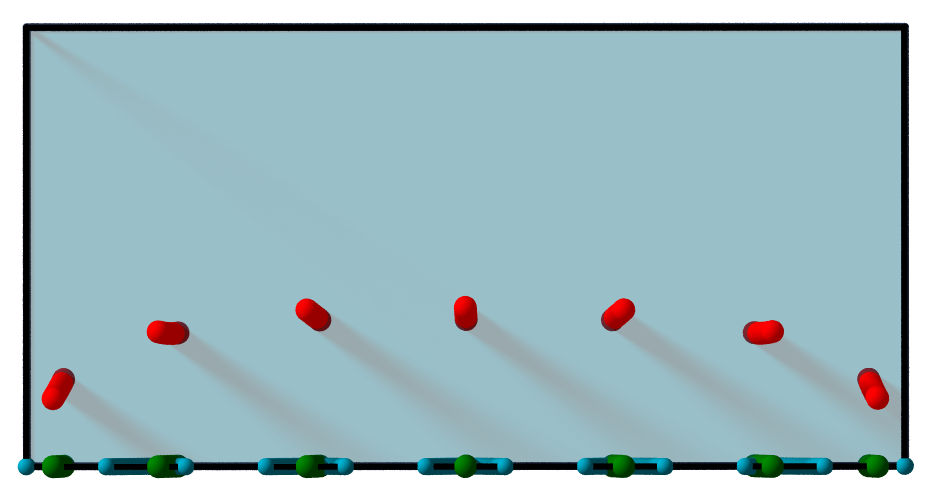}
  \put(24,24){\small$[\tau_3(t), \omega_3(t)]$}
\end{overpic}
\hfill \vrule width0pt\\[3ex]
\vrule width0pt\hfill
 \begin{overpic}[width=.99\textwidth,angle=0]{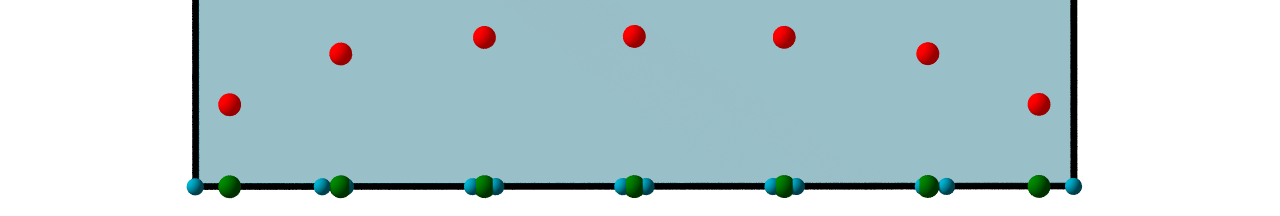}
 \put(1,15){(d)}
 \put(2,8){\fcolorbox{gray}{white}{\small$t=0.25$}}
\end{overpic}\hfill \vrule width0pt\\[1ex]
\vrule width0pt\hfill
 \begin{overpic}[width=.99\textwidth,angle=0]{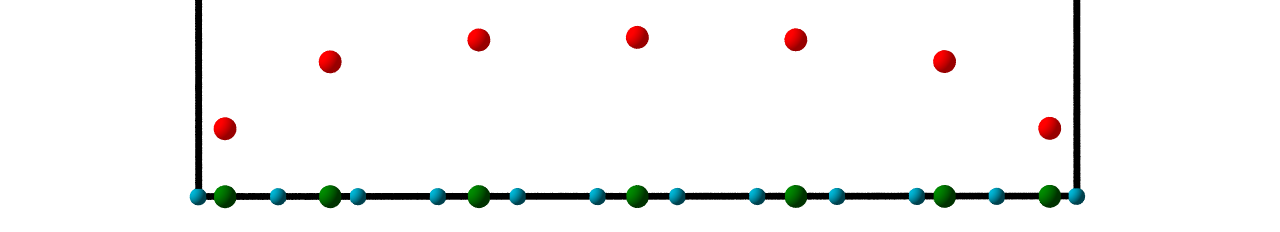}
 \put(1,15){(e)}
 \put(45,17){\small$[\tau_4,\omega_4]$}
 \put(2,8){\fcolorbox{gray}{white}{\small$t=1$}}
	\end{overpic}
\hfill \vrule width0pt\\[-4ex]
\Acaption{1ex}{Geodesic continuation of the Gaussian quadrature rules for $N=11$.
(a) The evolution of a Gaussian quadrature rule
from the source space $\St_{3,1}^n$ with $n=6$ is shown.
As the knot vector changes (blue paths), the nodes (green) and the weights (red) move continuously from 
their source position ($t=0$) to the target configuration $(t=1)$.
The evolution of $m=n+1=7$ Gaussian nodes and weights on $[a,b]$ is shown from the top view (b) and the side view (c).
(d) A zoom-in on the Gaussian quadrature rule for a specific time instant $t=0.25$.
The current nodes (green), weights (red) and knots (blue) are shown.
(e) The target uniform knot configuration and the target rule at $t=1$.}\label{fig:Continuation}
 \end{figure}

In this section, we show the results of our algorithm and
derive Gaussian quadrature rules for various target spaces $S_{3,2}^N$.
The evolution of the quadrature rules for $N=11$ is shown in Fig.~\ref{fig:Continuation}.
Starting with the Gaussian quadrature rule for a spline space with uniformly distributed double knots as a source rule,
the target rule for $S_{3,2}^N$ is traced by the homotopy continuation-based
Algorithm~\ref{algor}, deriving the optimal rules also for all the intermediate spaces.

\begin{figure}[!tb]
\vrule width0pt\hfill
 \begin{overpic}[width=.32\columnwidth,angle=0]{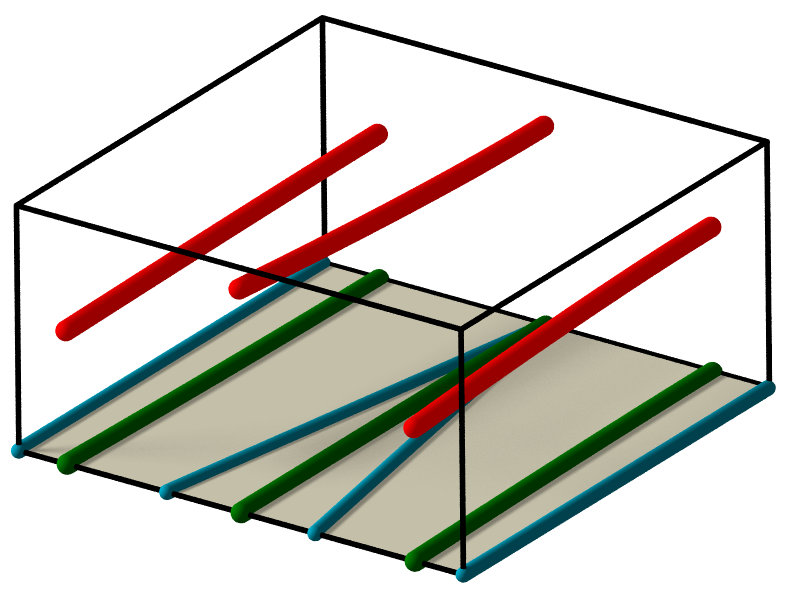}
   \put(10,50){$\omega_1(t)$}
   \put(0,71){\fcolorbox{gray}{white}{\small $N=3$}}
	\end{overpic}
 \hfill
 \begin{overpic}[width=.32\columnwidth,angle=0]{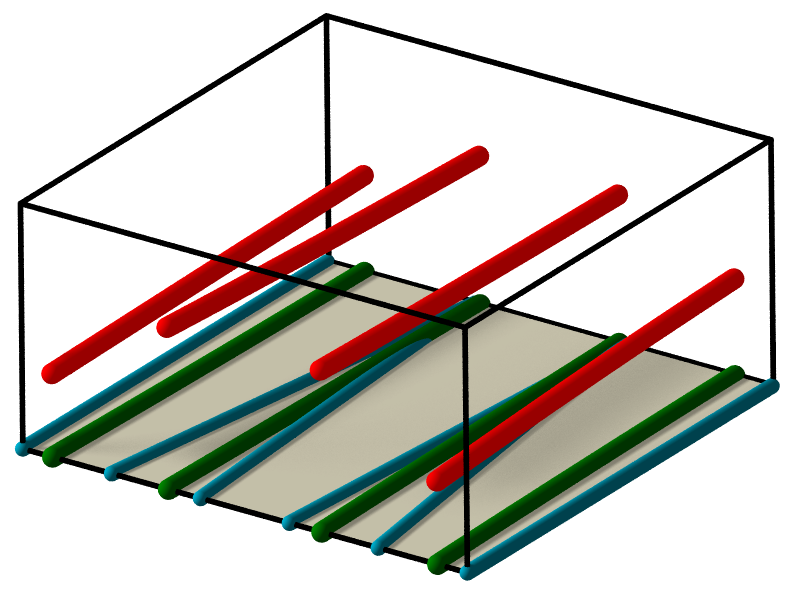}
 \put(0,71){\fcolorbox{gray}{white}{\small $N=5$}}
\end{overpic}
\hfill
 \begin{overpic}[width=.32\columnwidth,angle=0]{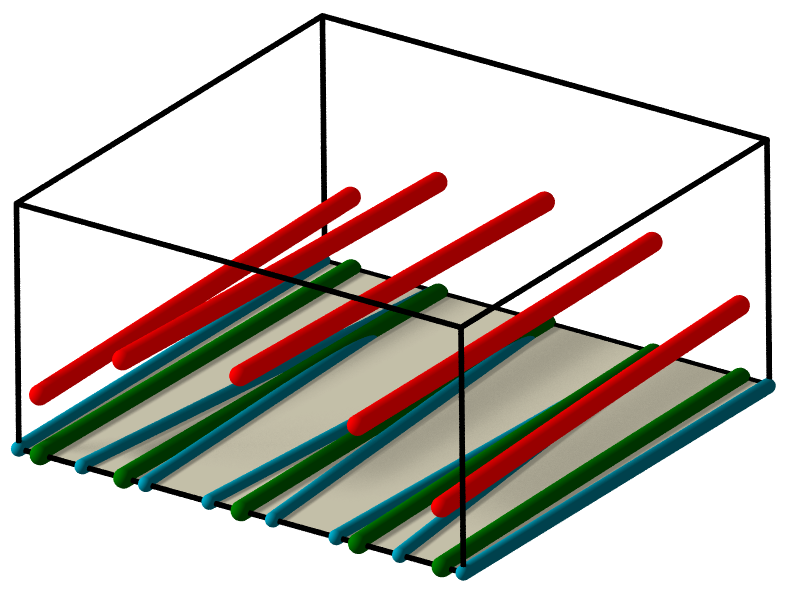}
  \put(0,71){\fcolorbox{gray}{white}{\small $N=7$}}
\end{overpic}\hfill \vrule width0pt\\[0.5ex]
\vrule width0pt\hfill
 \begin{overpic}[width=.32\columnwidth,angle=0]{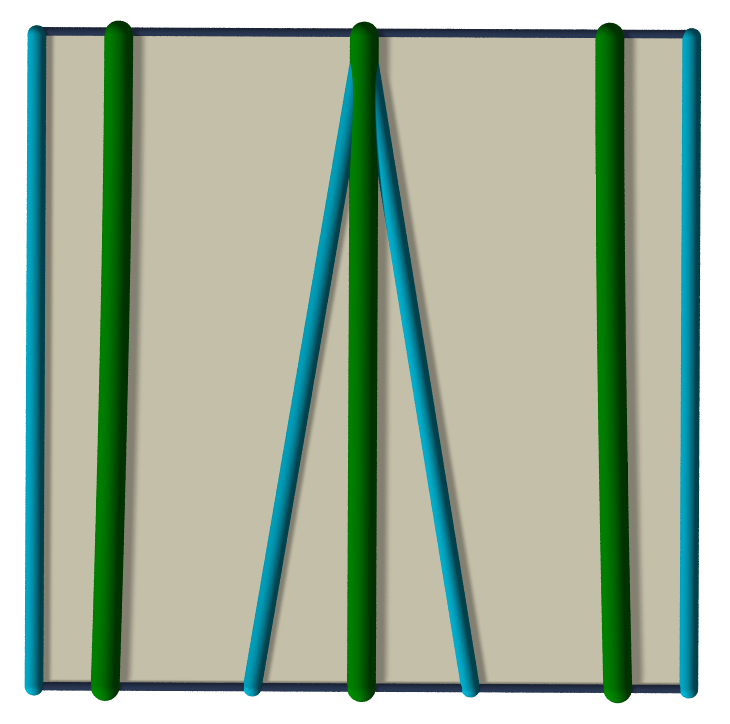}
  \put(13,99){\small $\tau_1(0)$}
 \put(11,-3){\small $\tau_1(1)$}
	\end{overpic}
 \hfill
 \begin{overpic}[width=.32\columnwidth,angle=0]{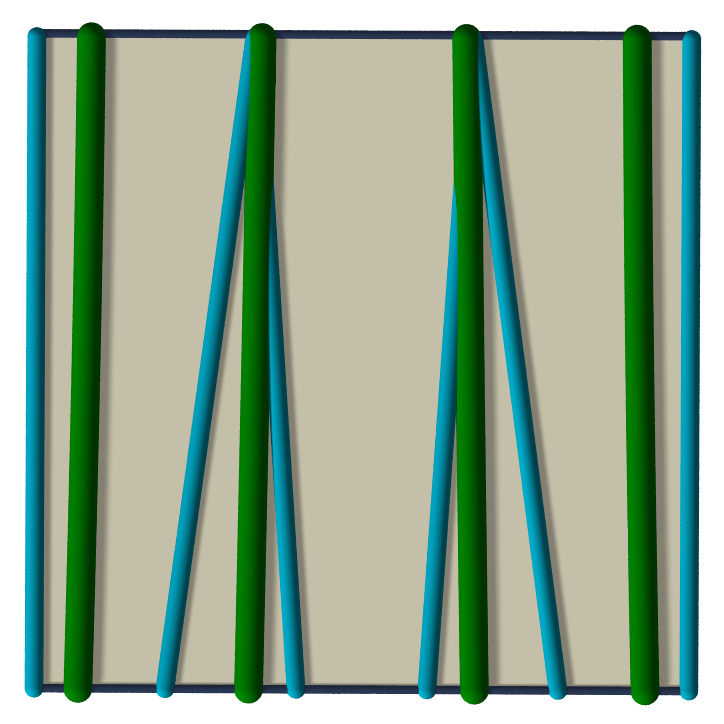}
\end{overpic}
\hfill
 \begin{overpic}[width=.32\columnwidth,angle=0]{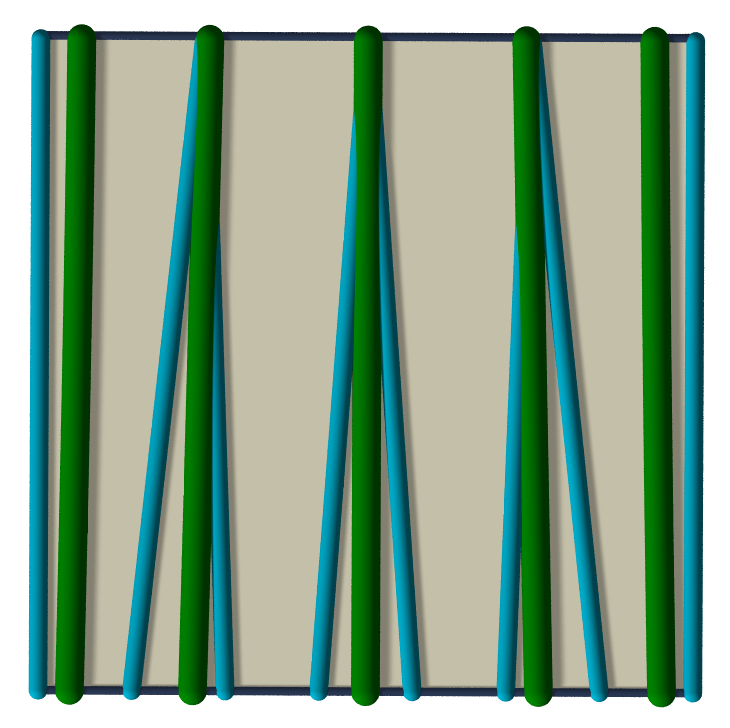}
\end{overpic}\hfill \vrule width0pt\\[1.5ex]
\vrule width0pt\hfill
 \begin{overpic}[width=.32\columnwidth,angle=0]{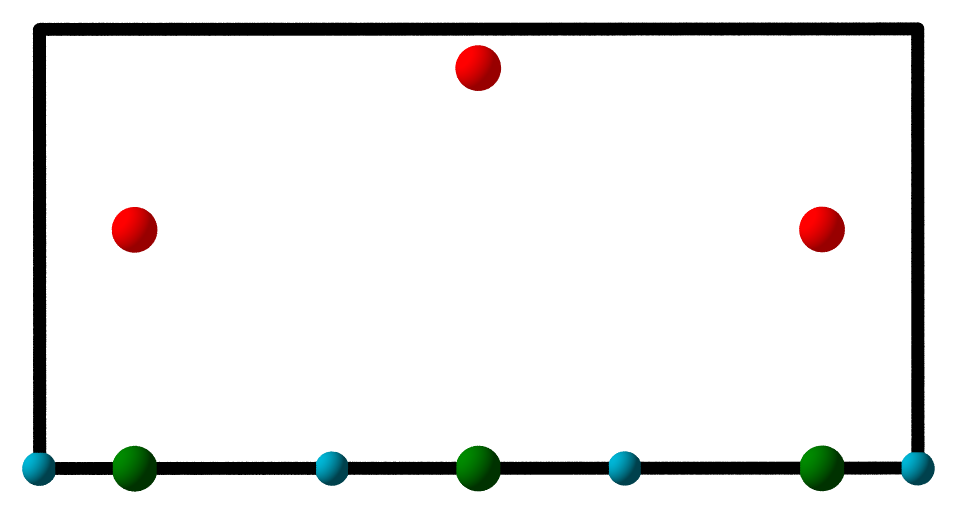}
 \put(7,35){$[\tau_1,\omega_1]$}
	\end{overpic}
 \hfill
 \begin{overpic}[width=.32\columnwidth,angle=0]{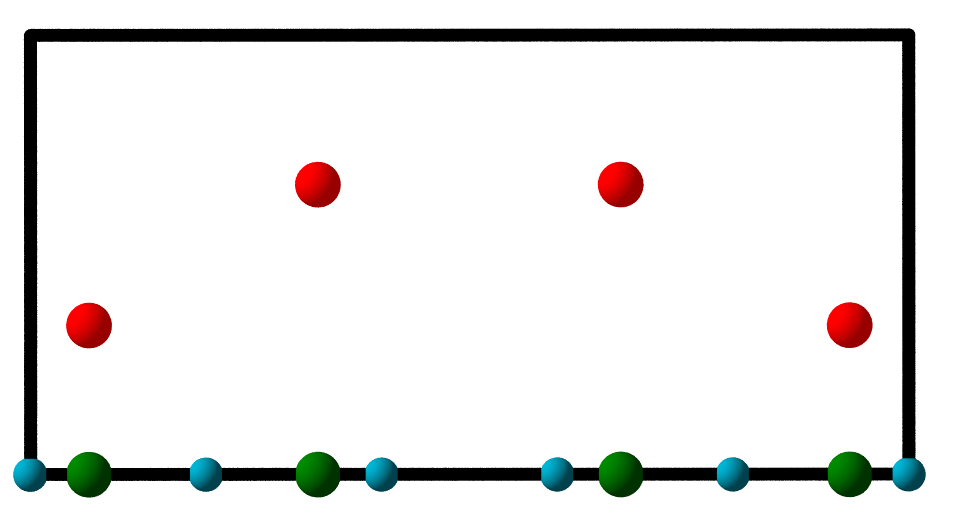}
\end{overpic}
\hfill
 \begin{overpic}[width=.32\columnwidth,angle=0]{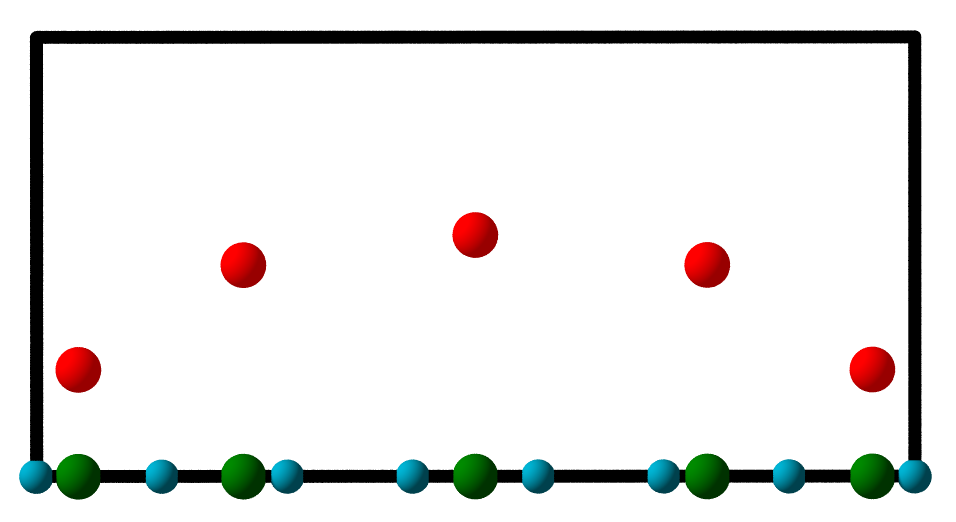}
\end{overpic}\hfill \vrule width0pt\\[4.ex]
\vrule width0pt\hfill
\begin{overpic}[width=.99\columnwidth,angle=0]{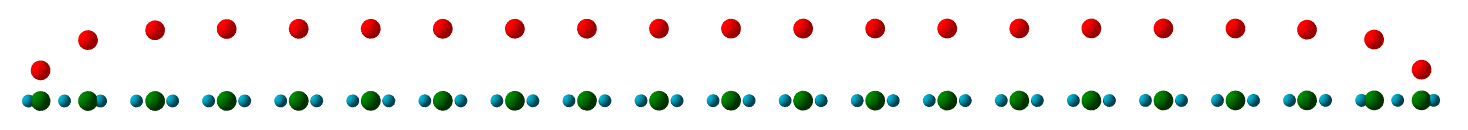}
\put(0,9){\fcolorbox{gray}{white}{\small $N=39$}}
	\end{overpic}
 \hfill \vrule width0pt\\[1ex]
 \vspace{-5ex}
\Acaption{1ex}{The evolution of Gaussian quadrature rules (\ref{quadrature}) for various source $N$ are shown.
Top: The evolution of the nodes (green), weights (red) for geodesic knot transformation (blue).
Middle: The evolution of nodes in time and the target rules for the time instant $t=1$.
Bottom: The rule for $N=39$; observe the convergence to the midpoint rule of Hughes et al.~\cite{Hughes-2010},
cf. Table~\ref{tabW}.}\label{fig:RuleVariousN}
 \end{figure}

  \begin{table}[!tb]
 \begin{center}
  \begin{minipage}{0.9\textwidth}
\caption{Nodes and weights for Gaussian quadrature rules (\ref{quadrature}) with double-precision for
uniform knot distribution for various target $N$ are shown.
To compare the results with algorithm of Nikolov \cite{Nikolov-2012}, the interval was set as $[a,b]=[0,1]$.
Due to the symmetry, only the first $[\frac{N+1}{4}]+1$ nodes and weights are displayed. The error $\|\br\|$
of the rule is measured as the Euclidean norm of the vector of the residues,
normalized by the system's dimension $N+3$, see (\ref{eq:Error}).}\label{tabW}
  \end{minipage}
\vspace{0.2cm}\\
\small{
\renewcommand{\arraystretch}{1.15}
\renewcommand{\tf}{\small}
\begin{tabular}{| c || l| l| l|}\hline
\rotatebox{0}{}
 $i$ &  $\tau_i$  & $\omega_i$ & $\|\br\|$\\\hline\hline
 & \multicolumn{2}{c|}{$N=3$} & \\
 1 & \tf 0.1086264370680297 & \tf 0.2720231005023455 & \multirow{2}{*}{$7.90^{-20}$}\\
 2 & \tf 0.5                & \tf 0.4559537989953090 & \\\hline\hline
 & \multicolumn{2}{c|}{$N=5$} & \\
 1 & \tf 0.0669578918742195 & \tf 0.1698605936669416 & \multirow{2}{*}{$1.04^{-19}$}\\
 2 & \tf 0.3275898516368645 & \tf 0.3301394063330584 & \\\hline\hline
  & \multicolumn{2}{c|}{$N=7$} & \\
 1 & \tf 0.0479188107803577 & \tf 0.1216810800700958 & \multirow{3}{*}{$1.95^{-18}$}\\
 2 & \tf 0.2358921494969001 & \tf 0.2408185184939348 & \\
 3 & \tf 0.5                & \tf 0.2750008028719389 & \\\hline\hline
   & \multicolumn{2}{c|}{$N=9$} & \\
 1 & \tf 0.0372757529111283 & \tf 0.0946622477445919 & \multirow{3}{*}{$2.08^{-18}$}\\
 2 & \tf 0.1835904624135774 & \tf 0.1876252194189693 & \\
 3 & \tf 0.3904233866079767 & \tf 0.2177125328364388 & \\\hline\hline
    & \multicolumn{2}{c|}{$N=11$} & \\
 1 & \tf 0.0304987043023585 & \tf 0.0774523185174377 & \multirow{4}{*}{$6.68^{-18}$}\\
 2 & \tf 0.1502181009517147 & \tf 0.1535325192913209 & \\
 3 & \tf 0.3195393932155687 & \tf 0.1783894870783702 & \\
 4 & \tf 0.5                & \tf 0.1812513502257421 & \\\hline
  \multicolumn{4}{c}{$\vdots$}\\\hline
    & \multicolumn{2}{c|}{$N=39$} &\\
 1 & \tf 0.0086022074347388 & \tf 0.0218455595269063 & \multirow{11}{*}{$1.02^{-17}$}\\
 2 & \tf 0.0423693959303822 & \tf 0.0433045545577068 & \\
 3 & \tf 0.0901289847662636 & \tf 0.0503213631747089 & \\
 4 & \tf 0.1410569521267253 & \tf 0.0512021143533085 & \\
 5 & \tf 0.1923101843694322 & \tf 0.0512756766459810 & \\
 6 & \tf 0.2435899416018961 & \tf 0.0512815446928528 & \\
 7 & \tf 0.2948718106031808 & \tf 0.0512820110347811 & \\
 8 & \tf 0.3461538474036372 & \tf 0.0512820480845737 & \\
 9 & \tf 0.3974358975351839 & \tf 0.0512820510280155 & \\
10 & \tf 0.4487179487257872 & \tf 0.0512820512617426 & \\
11 & \tf 0.5                & \tf 0.0512820512788446 & \\\hline
\end{tabular}
}
\end{center}
\end{table}

 \begin{figure}[!tbh]
\vrule width0pt\hfill
 \begin{overpic}[width=.33\columnwidth,angle=0]{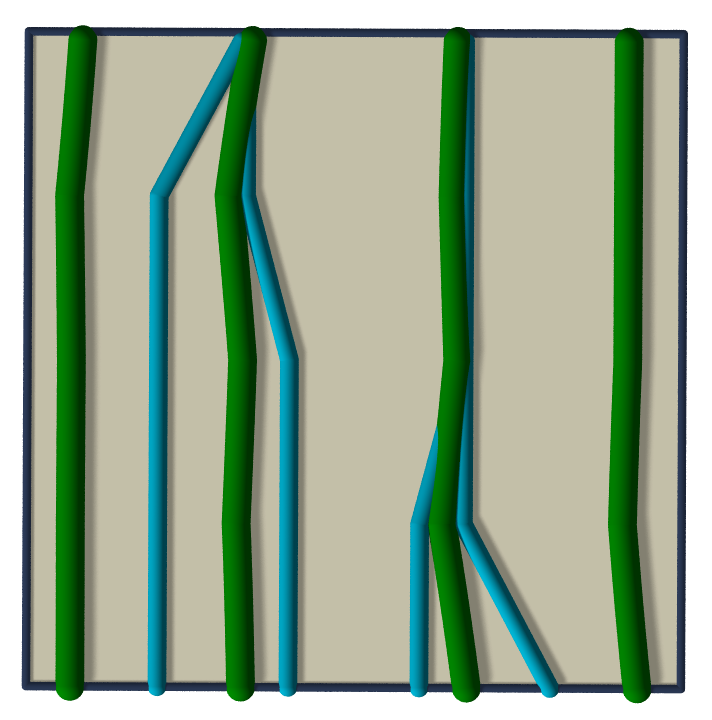}
   \put(30,-6){$(1,2,3,4)$}
	\end{overpic}
 \hfill
 \begin{overpic}[width=.33\columnwidth,angle=0]{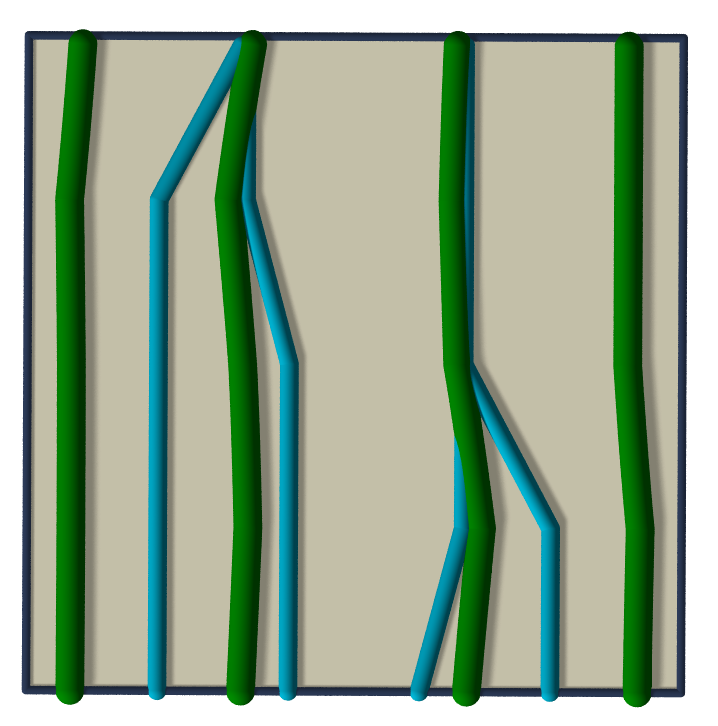}
    \put(30,-6){$(1,2,4,3)$}
\end{overpic}
\hfill
 \begin{overpic}[width=.33\columnwidth,angle=0]{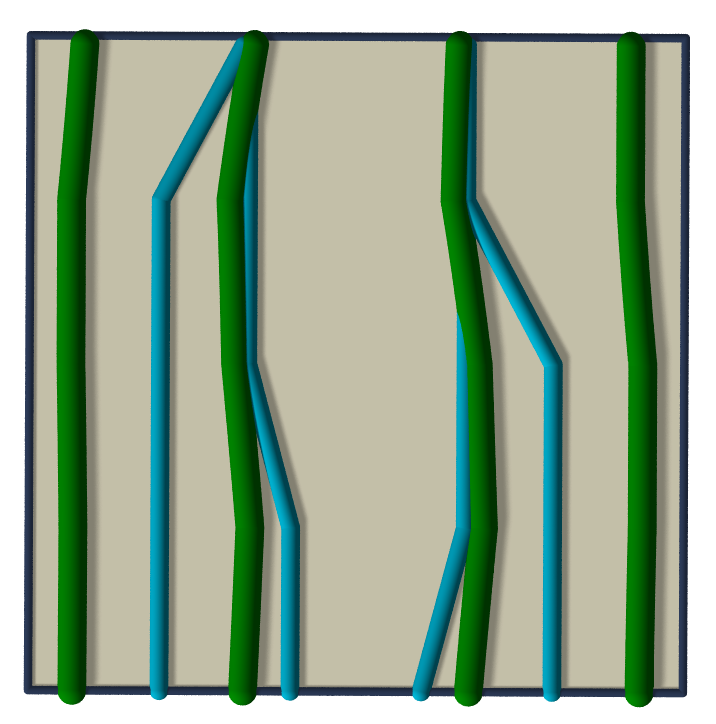}
   \put(30,-6){$(1,4,2,3)$}
   \put(100,95){\vector(0,-1){90}}
    \put(93,100){\small$t=0$}
    \put(93,-6){\small$t=1$}
\end{overpic}\hfill \vrule width0pt\\[3ex]
\vrule width0pt
 \begin{overpic}[width=.33\columnwidth,angle=0]{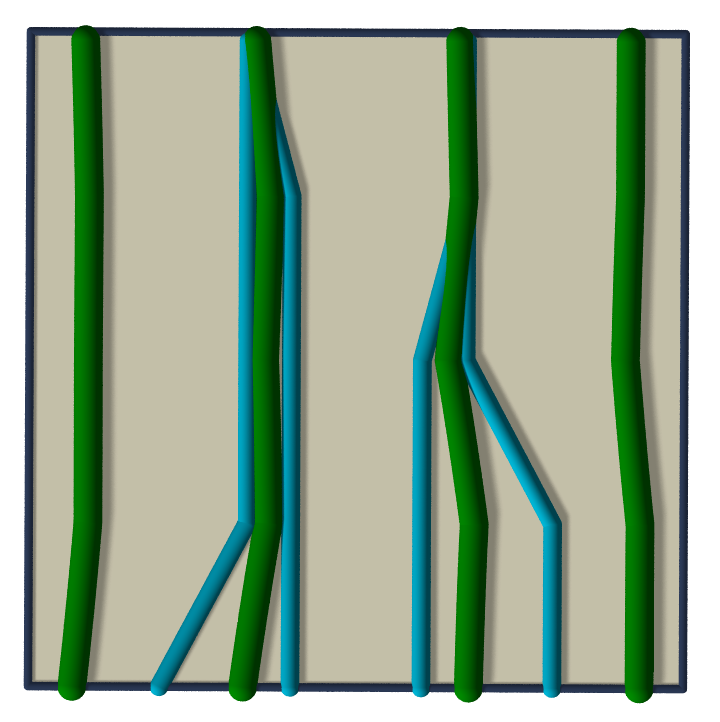}
        \put(30,-6){$(2,3,4,1)$}
	\end{overpic}
 \hfill
 \vrule width0pt\\[3ex]
\vrule width0pt
 \begin{overpic}[width=.33\columnwidth,angle=0]{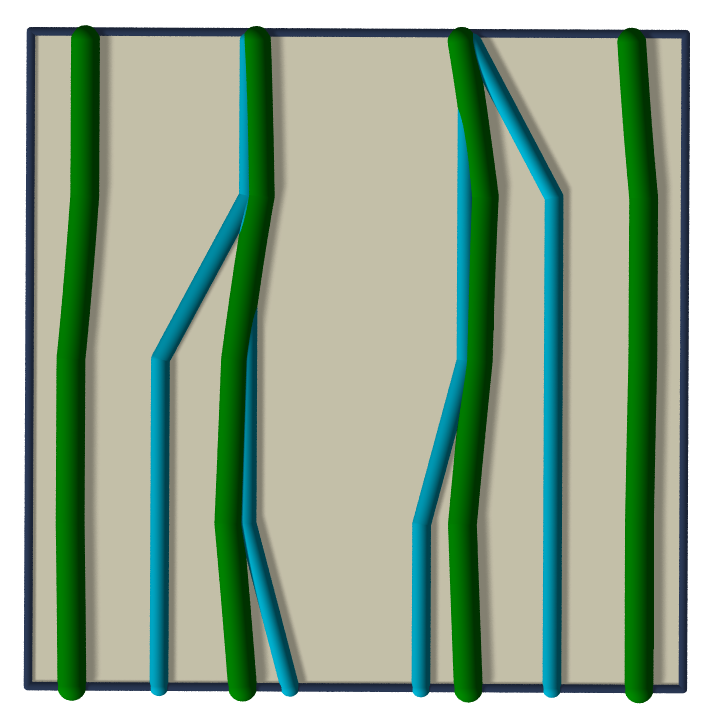}
 \put(30,-6){$(4,1,3,2)$}
 \raisebox{0.30\columnwidth}{
   \put(120,0){\includegraphics[width=.55\columnwidth,angle=0]{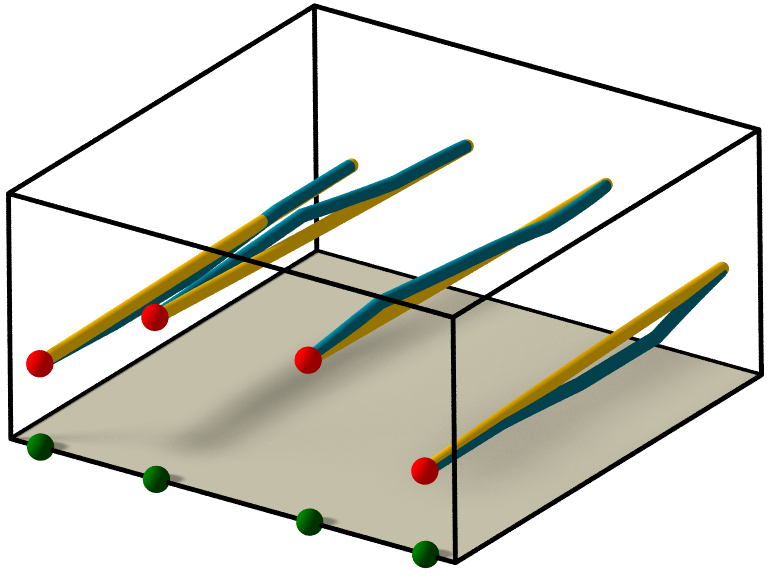}}
   \put(105,1.5){$(4,1,3,2)$ vs. geodesic}
   \put(170,7){\vector(0,4){57}}
   \put(170,7){\vector(-1,2){24.3}}
   \put(225,0){$t=1$}
   \put(280,34){$t=0$}
   }
   \put(100,40){
   \begin{minipage}{0.65\textwidth}
    \Acaption{1ex}{Path independence. The evolution of the Gaussian quadrature rules (\ref{quadrature}) for various modifications
    of the knot vector for $N=5$ is shown. At every time instant, only one knot changes its position.
    Five different evolutions of the Gaussian nodes are shown (green). Starting with two double knots, the quadruplet encodes the particular order of
    the four knots to be moved towards their target position (uniform spacing at $t=1$). The 3D figure compares the $(4,1,3,2)$-evolution of the weights (blue)
    with the geodesic evolution (yellow), cf. Fig.~\ref{fig:RuleVariousN}. All paths derive the same target rule ($t=1$) visualized by red dots (weights)
    and green dots (nodes).}\label{fig:AlongEdge}
   \end{minipage}
   }
\end{overpic}
\hfill \vrule width0pt\\[1ex]
 \end{figure}

 \begin{table}[!tb]
 \begin{center}
  \begin{minipage}{0.9\textwidth}
\caption{Path independence. The target quadrature rule obtained by different along-the-edge
paths of the knot vector. The results for five paths shown in Fig.~\ref{fig:AlongEdge} are displayed
with twenty digits accuracy. The target nodes and weights match up to eighteen digits.}\label{tab2}
  \end{minipage}
\vspace{0.2cm}\\
\small{
\renewcommand{\arraystretch}{1.2}
\renewcommand{\tf}{\small}
\begin{tabular}{| c || l| l|}\hline
\rotatebox{0}{}
  &  $\tau_1$, $\tau_2$  & $\omega_1$, $\omega_2$ \\\hline\hline
 \multirow{2}{*}{$(1,2,3,4)$} & \tf 0.06695789187421950918 & \tf 0.16986059366694164265 \\
                              & \tf 0.32758985163686446374 & \tf 0.33013940633305835725 \\\hline\hline
 \multirow{2}{*}{$(1,2,4,3)$} & \tf 0.06695789187421950918 & \tf 0.16986059366694164265 \\
                              & \tf 0.32758985163686446374 & \tf 0.33013940633305835725 \\\hline\hline
 \multirow{2}{*}{$(1,4,2,3)$} & \tf 0.06695789187421950918 & \tf 0.16986059366694164265 \\
                              & \tf 0.32758985163686446374 & \tf 0.33013940633305835725 \\\hline\hline
 \multirow{2}{*}{$(2,3,4,1)$} & \tf 0.06695789187421950915 & \tf 0.16986059366694164265 \\
                              & \tf 0.32758985163686446341 & \tf 0.33013940633305835787 \\\hline\hline
 \multirow{2}{*}{$(4,1,3,2)$} & \tf 0.06695789187421950915 & \tf 0.16986059366694164265 \\
                              & \tf 0.32758985163686446342 & \tf 0.33013940633305835721 \\\hline
\end{tabular}
}
\end{center}
\end{table}

 The evolution of the nodes is visualized as
a set of planar curves, the evolution of weights is represented by the corresponding space trajectories.
The error of the rule $\Q$ is measured in terms of the Euclidean norm of the vector of the residues of the
system (\ref{eq:IniSystem}), normalized by the dimension of the system
\begin{equation}\label{eq:Error}
\|\br\| = \frac{1}{N+3}(\sum_{i=1}^{N+3} (\Q_a^b[D_i] - I[D_i])^2)^{\frac{1}{2}}.
\end{equation}
Our results coincide with those of \cite{Nikolov-2012} (up to twelve digits precision shown therein),
where the correct layout of nodes was conjectured as
\begin{equation}\label{NikolovNP}
\bp = (1,1,0,1,0,1 \dots, 1,0,1,0,1,1),
\end{equation}
 that is, the first two boundary elements contain a single node each, whilst the middle elements
 follow a node-gap pattern. Our algorithm yielded the same nodal pattern, see Fig.~\ref{fig:RuleVariousN}.
 Their algorithm is a trial-error scheme that chooses the first node and, under the assumption of the nodal pattern (\ref{NikolovNP}),
 computes iteratively the remaining nodes and weights. Based on the error of the rule,
 the first node is modified and the next iteration is invoked. Such an approach relies heavily on the initial node
 and, since the problem is highly non-linear. In our experience, 
 there is no guarantee that the scheme of \cite{Nikolov-2012} will converge to a correct rule.
 We emphasize that our approach derives the layout of nodes (\ref{NikolovNP}) automatically and
 is general because it is not limited to only uniform target knot sequences.

  Also observe the convergence to the midpoint rule of Hughes et al.~\cite{Hughes-2010},
 see Table~\ref{tabW}, where the limit weight for $N=39$ ($n = 20$) is
 \begin{equation}\label{eq:limits}
\frac{2}{39} = 0.\overline{051282}
\end{equation}
The data in Table~\ref{tabW} are shown with double-precision (16 decimal digits).
Observe the slow convergence, e.g., for $N=39$, the
results meet the half-point-rule limit with only five digits of accuracy.
This fact is not a limitation of our work, on the contrary, it shows that the midpoint rule is an approximation
of our exact rule and can be used only when $N$ is large enough.

\begin{figure}[!tb]
\vrule width0pt\hfill
 \begin{overpic}[width=.32\columnwidth,angle=0]{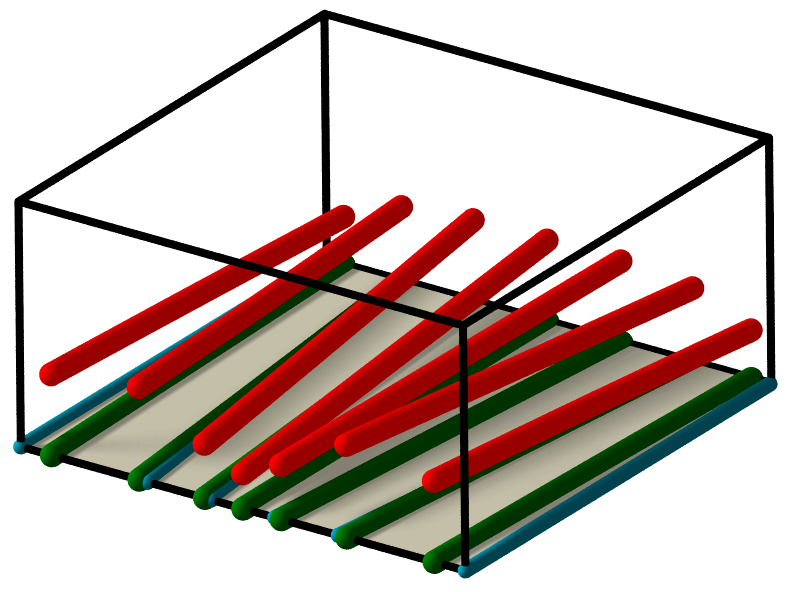}
   \put(0,71){\fcolorbox{gray}{white}{\small $q=2$}}
	\end{overpic}
 \hfill
 \begin{overpic}[width=.32\columnwidth,angle=0]{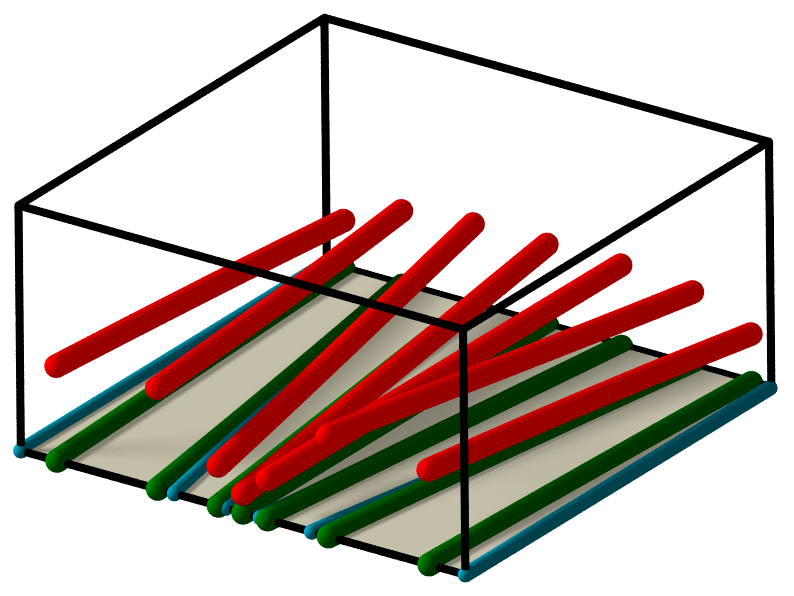}
 \put(0,71){\fcolorbox{gray}{white}{\small $q=3$}}
\end{overpic}
\hfill
 \begin{overpic}[width=.32\columnwidth,angle=0]{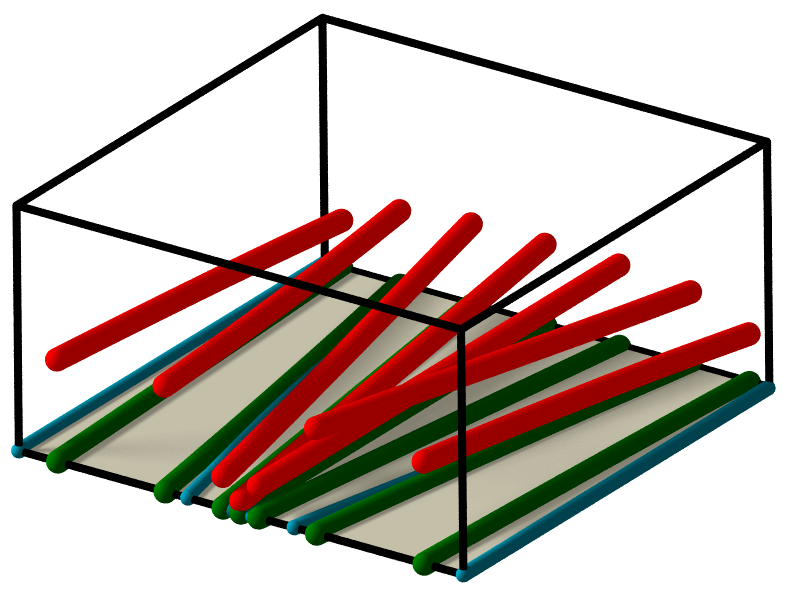}
 \put(0,71){\fcolorbox{gray}{white}{\small $q=4$}}
\end{overpic}\hfill \vrule width0pt\\[0.5ex]
\vrule width0pt\hfill
 \begin{overpic}[width=.32\columnwidth,angle=0]{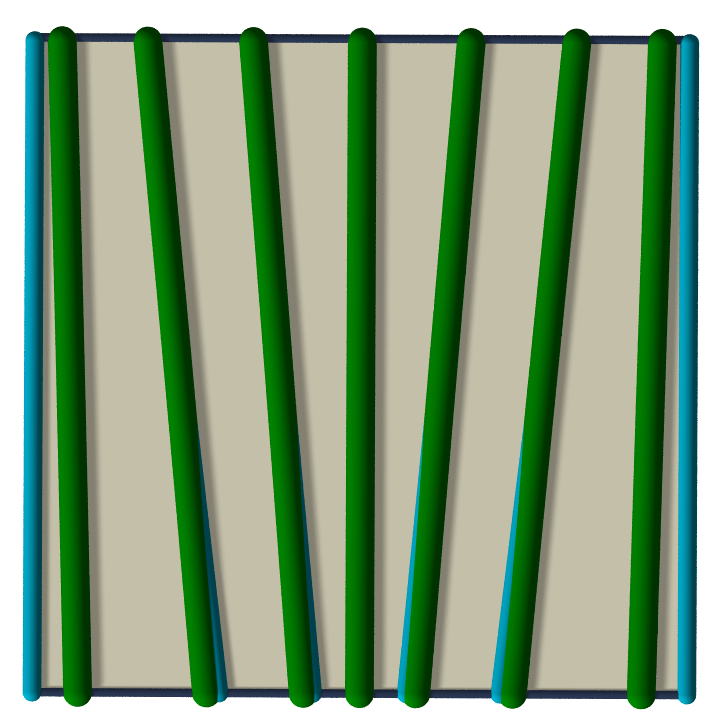}
 \put(13,99){\small $\tau_2(0)$}
 \put(19,-4){\small $\tau_2(1)$}
	\end{overpic}
 \hfill
 \begin{overpic}[width=.32\columnwidth,angle=0]{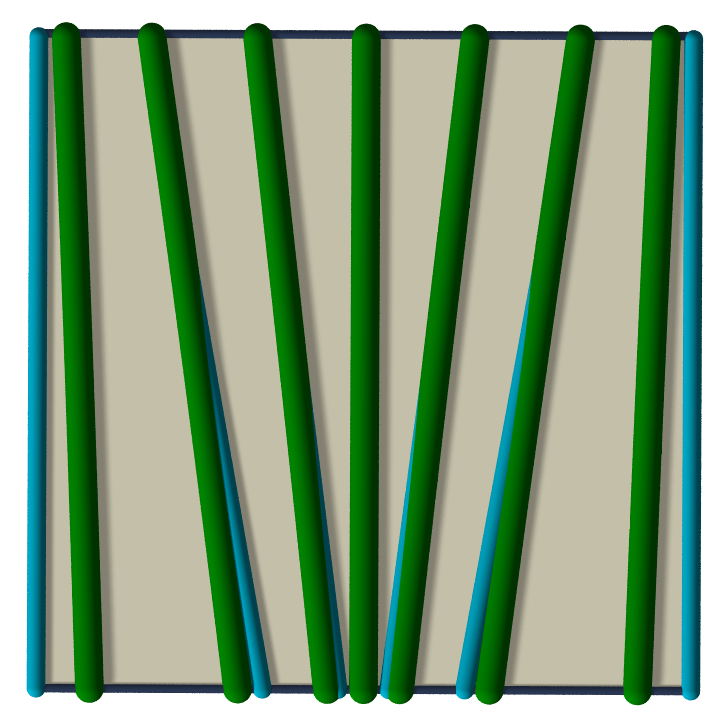}
\end{overpic}
\hfill
 \begin{overpic}[width=.32\columnwidth,angle=0]{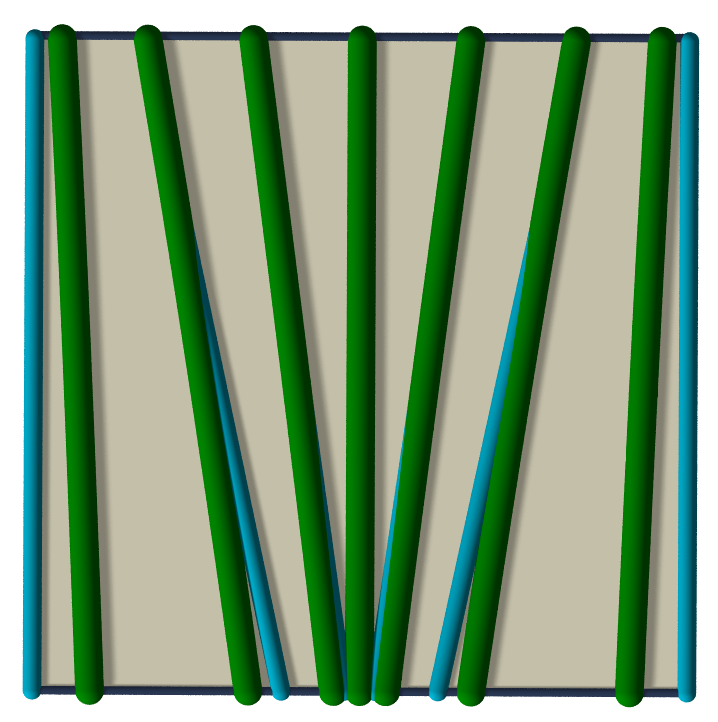}
\end{overpic}\hfill \vrule width0pt\\[1.5ex]
\vrule width0pt\hfill
 \begin{overpic}[width=.32\columnwidth,angle=0]{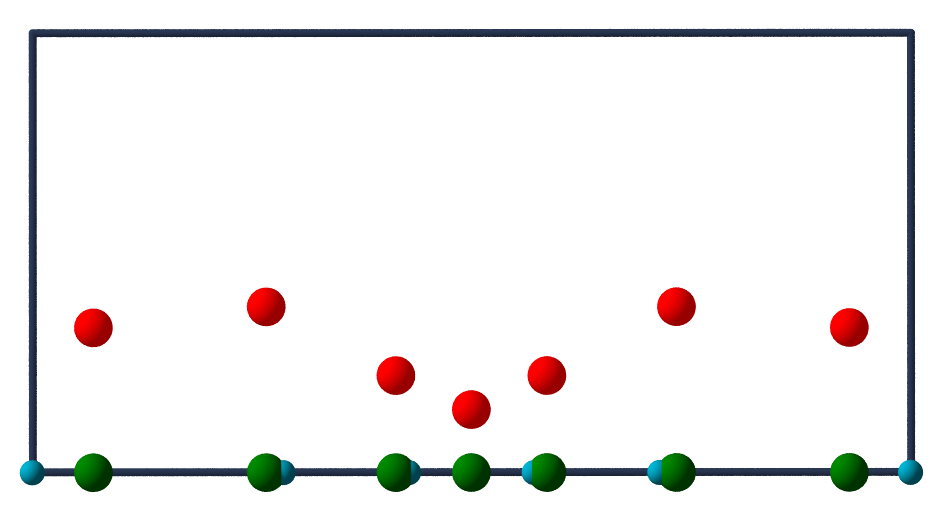}
	\end{overpic}
 \hfill
 \begin{overpic}[width=.32\columnwidth,angle=0]{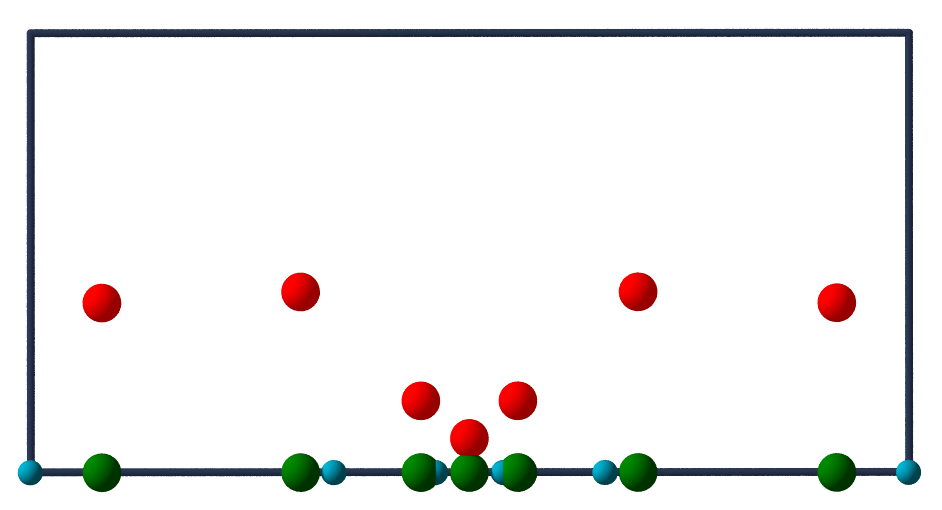}
\end{overpic}
\hfill
 \begin{overpic}[width=.32\columnwidth,angle=0]{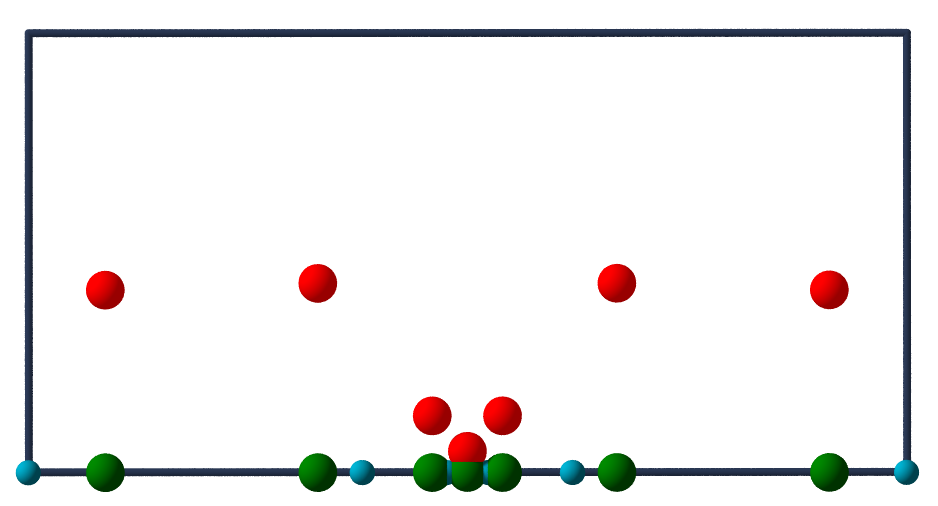}
\end{overpic}\hfill \vrule width0pt\\[-1.ex]
\Acaption{1ex}{The evolution of the Gaussian rules for $C^1$ spaces with
non-uniform target knot sequences.
Five interior source double knots (blue) are moved towards the interval's center. The parameter $q$ is
the shrinking ratio; it is the ratio of lengths of two neighboring elements when moving towards the center.
Top: The evolution of the rule in time visualized in 3D: the trajectories of the weights (red), nodes (green) and knots (blue)
are shown.
Middle: The evolution of the nodes. For all $q$ shown, $\tau_2$ starts ($t=0$) in the second elements but
ends ($t=1$) in the first element. Bottom: The target Gaussian rules for particular non-uniform knot sequences.
}\label{fig:Shrinking}
 \end{figure}

The evolution of the Gaussian nodes and weights for various $N$ is shown in Fig.~\ref{fig:RuleVariousN}.
The modification of the knots is geodesic, see Section~\ref{sec:Splines}, i.e., every knot transforms
to its target destination along the shortest possible path.
As previously discussed in Section~\ref{sec:HomotopyQuad}, the source and the target rules are both unique.
Therefore, any path deforming the uniform double knot sequence $\XXt_n$ (\ref{eq:XXt})
into the uniform knot sequence $\XX_N$ (\ref{eq:XX}) must produce the same rule. The results of
five random along-the-edge knot modifications are shown in Fig.~\ref{fig:AlongEdge}.
Contrary to the geodesic knot deformation,
only one knot moves while the other stay fixed during the evolution.
Depending on the order of how the knots are moved, there exist $(N+5)!$ possible paths.
The results of the target rule for various along-the-edge knot vector deformations are shown in Table~\ref{tab2}.
The results are independent on the path which validates the correctness of our algorithm.


\begin{figure}[!tb]
\vrule width0pt\hfill
 \begin{overpic}[height=.49\columnwidth,angle=0]{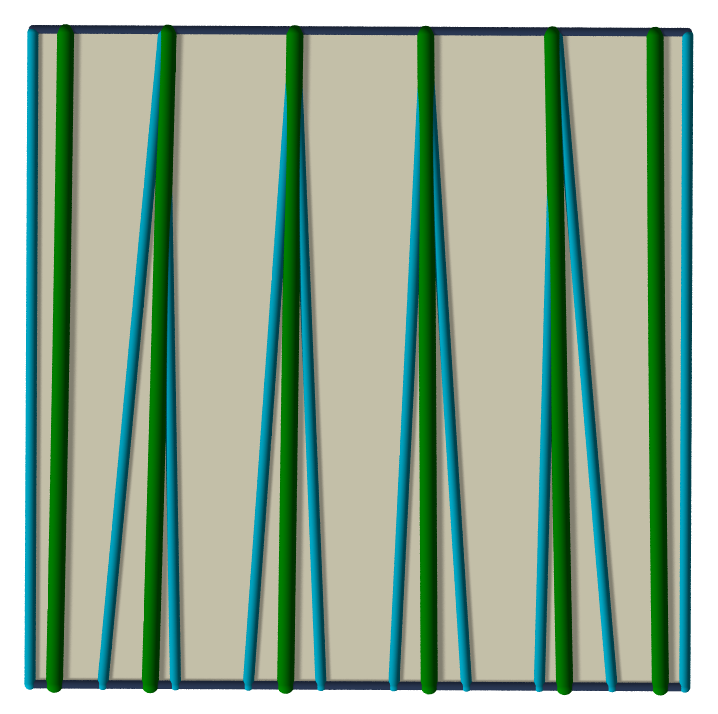}
   \put(30,105){\fcolorbox{gray}{white}{\small $\St_{3,1}^n$ with $n=5$}}
   \put(80,-8){\fcolorbox{gray}{white}{\small $S_{3,2}^N$ with $N=9$}}
	\end{overpic}
\hfill
 \begin{overpic}[height=.49\columnwidth,angle=0]{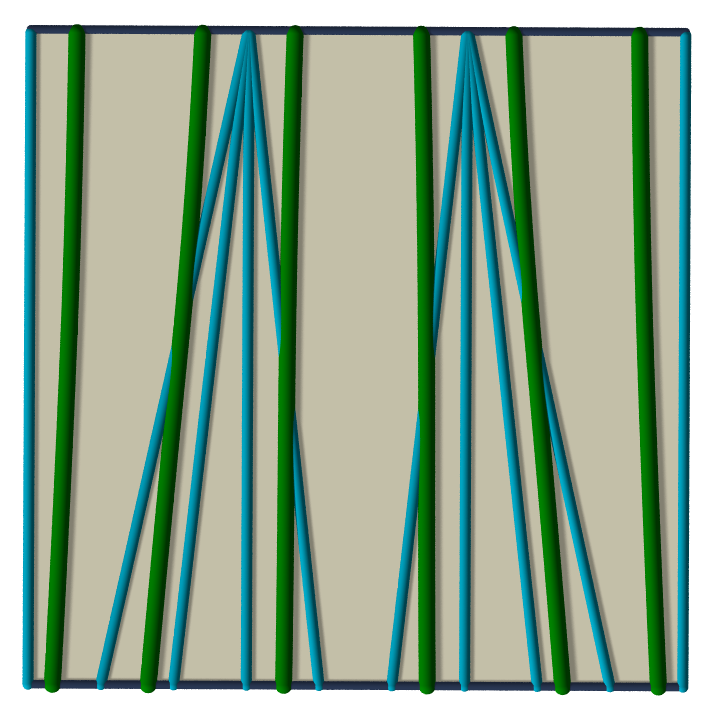}
 \put(30,105){\fcolorbox{gray}{white}{\small $\St_{3,-1}^n$ with $n=3$}}
\end{overpic}\hfill \vrule width0pt\\[5ex]
\vrule width0pt\hfill
 \begin{overpic}[width=.99\columnwidth,angle=0]{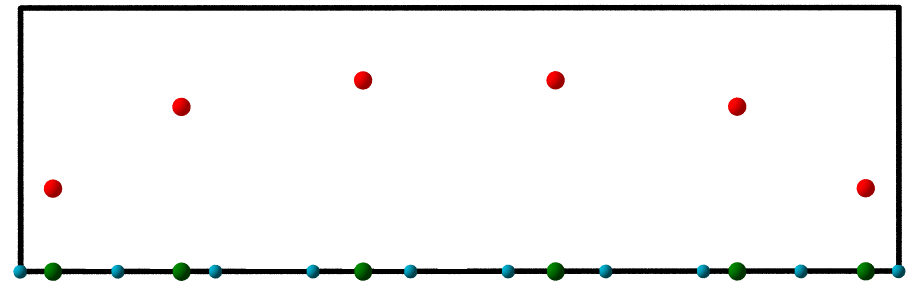}
 \put(40,15){\fcolorbox{gray}{white}{\small target rule $\Q(1)$}}
 \put(20,23){$[\tau_2,\omega_2]$}
\end{overpic}
\hfill \vrule width0pt\\[-4ex]
\Acaption{1ex}{Evolution of the Gaussian rule from two different source rules is shown:
the Gaussian rule for $\St_{3,1}^n$ with $n=5$ uniform elements (left)
and the classical Gaussian rule for polynomials with $n=3$ elements (right).
The target space is $S_{3,2}^N$, the space of $C^2$ cubics with $N=9$ uniform elements.
The same target rule (bottom) was derived by Algorithm~\ref{algor} with the accuracy of $\|\br\|<10^{-16}$,
see (\ref{eq:Error}).}\label{fig:Gauss}
 \end{figure}

An example of Gaussian rules for non-uniform target knot sequences is shown in Fig.~\ref{fig:Shrinking}.
The space is kept $C^1$, i.e., all the knots keep their multiplicity two. The target knot sequence was required to have
larger elements close to the boundary and smaller ones in the interval's center. Unlike our previous result \cite{Quadrature31-2014}
where finer spacing close to the boundary guaranteed the same nodal pattern for the whole family of
symmetrically stretched knot sequences, here we can observe a change of the nodal pattern of the rule.

Throughout the paper, we used the rule of Nikolov \cite{Nikolov-1996} as our source rule.
However, one may start, e.g., with the classical Gaussian rule for polynomials, see Fig.~\ref{fig:Gauss}.
In such a case, the knots have multiplicity four and one needs to count the proper number of initial elements
to get the desired number of target elements $N$. This counting obeys the rule of Micchelli \cite{Micchelli-1977}, see (\ref{eq:Micchelli}).
The results show the stability of our Algorithm~\ref{algor}. The results from Fig.~\ref{fig:Gauss} both correspond to values shown in
Table~\ref{tabW} for $N=9$. Both target rules derived from different source rules match up to sixteen decimal digits.

\section{Conclusion}\label{sec:conl}

We have introduced a new methodology to compute Gaussian quadrature rules.
Starting with a known Gaussian quadrature rule for a specific spline space, the rule is interpreted
as a point in a higher dimensional space, a zero of a specific polynomial system.
A homotopy continuation-based algorithm has been presented that
numerically traces the Gaussian quadrature rule as the knot vector,
and consequently the whole spline space, is continuously modified.
We have recovered the Gaussian quadrature rule of Nikolov~\cite{Nikolov-2012}
for the uniform $C^2$ cubic spline spaces, a rule that was derived under an assumption
of the conjectured nodal pattern (\ref{NikolovNP}).

The continuation approach shows the connection between different Gaussian rules.
This concept eliminates the need for the construction of particular rules for different kinds
of non-uniform knot vectors. Thus, work such as \cite{Quadrature31-2014,Quadrature51-2014}
is hereafter needed only to validate numerical results obtained with this method.
We plan to further exploit the proposed technique to derive
the Gaussian rules for the particular spline spaces of a high interest such as $S_{4,0}$ with
non-uniform knot sequences, or $S_{6,1}$, spaces frequently appearing in the Galerkin method
when building the mass and stiffness matrices.


\bibliographystyle{plain}
\bibliography{HomotopyQuadratureArXiv}



\end{document}